\documentclass{ps}
\usepackage{mathrsfs}
\usepackage{amsthm}
\usepackage{amsmath}
\usepackage{amssymb}
\usepackage{esint}
\usepackage[numbers]{natbib}
\usepackage{bbm}

\theoremstyle{plain}
\newtheorem{assumption}{Assumption}

\renewcommand{\hat}{\widehat}
\renewcommand{\tilde}{\widetilde}
\renewcommand{\bar}{\overline}

\allowdisplaybreaks[3]

\makeatother

\global\long\def\phi{\varphi}
\global\long\def\epsilon{\varepsilon}
\global\long\def\theta{\vartheta}
\global\long\def\E{\mathbb{E}}
\global\long\def\N{\mathbb{N}}
\global\long\def\Z{\mathbb{Z}}
\global\long\def\R{\mathbb{R}}
\global\long\def\F{\mathcal{F}}
\global\long\def\le{\leqslant}
\global\long\def\ge{\geqslant}
\global\long\def\1{\mathbbm 1}
\global\long\def\d{\mathrm{d}}

\global\long\def\Cov{\operatorname{Cov}}
\global\long\def\Span{\operatorname{span}}
\global\long\def\Id{\operatorname{Id}}
\global\long\def\subset{\subseteq}

\global\long\def\bull{{\scriptstyle \bullet}}
\global\long\def\supp{\operatorname{supp}}
\global\long\def\sgn{\operatorname{sign}}
\begin{document}
\title{Adaptive confidence bands for Markov chains and diffusions:
 Estimating the invariant measure and the drift}\thanks{The authors acknowledge intensive and very 
helpful discussions with Richard Nickl. 
J.S. thanks the European Research Council (ERC) for support under Grant No. 647812.
M.T. is grateful to the Statistical Laboratory of the University of Cambridge for its hospitality during a visit 
from February to March 2014, where this research was initiated, and to the Deutsche Forschungsgemeinschaft (DFG) for the research fellowship TR 1349/1-1. Part of the paper was carried out while M.T. was employed at
the Humboldt-Universit\"at zu Berlin.}
\runningtitle{Adaptive confidence bands for Markov chains and diffusions}

\author{Jakob S\"ohl}\address{Statistical Laboratory, Department of Pure Mathematics and Mathematical Statistics, University of Cambridge, CB3 0WB Cambridge, UK. Email: j.soehl@statslab.cam.ac.uk}
\author{Mathias Trabs}\address{Department of Mathematics, University of Hamburg, Bundesstra\ss e 55, 20146 Hamburg, Germany. Email: 
mathias.trabs@uni-hamburg.de}
%
%

\begin{abstract} 
  As a starting point we prove a functional central limit theorem for estimators of the invariant measure of a geometrically ergodic Harris-recurrent Markov chain in a multi-scale space. This allows to construct confidence bands for the invariant density with optimal (up to undersmoothing) $L^{\infty}$-diameter by using wavelet projection estimators. In addition our setting applies to the drift estimation of diffusions observed discretely with fixed observation distance. We prove a functional central limit theorem for estimators of the drift function and finally construct adaptive confidence bands for the drift by using a completely data-driven estimator.
\end{abstract}
%
%
\maketitle

\section*{Introduction}

Diffusion processes are prototypical examples of the theory of stochastic
differential equations as well as of continuous time Markov processes. At
the same time diffusions are widely used in applications, for instance, to model
molecular movements, climate data or in econometrics. Focusing on Langevin
diffusions, we will consider the solution of the stochastic differential
equation
\begin{equation*}
  \d X_t=b(X_t)\d t+\sigma\d W_t,\quad t\ge0,
\end{equation*}
with unknown drift function $b\colon\R\to\R$, a volatility parameter $\sigma>0$ and with a Brownian motion $W=\{W_t:t\ge0\}$. The problem of statistical estimation based on discrete observations from this model is embedded into the framework of geometrically ergodic Harris-recurrent Markov chains. We study the estimation of the invariant density of such Markov chains. The drift function $b$ depends nonlinearly on the invariant density $\mu$ so that the two estimation problems of $b$ and $\mu$ are closely related. We prove functional central limit theorems for estimators of both $b$ and $\mu$ in multi-scale spaces. This allows the construction of confidence bands for $\mu$. Owing to the nonlinear dependence the construction of confidence bands for $b$ is more involved. In this more difficult situation and by using a self-similarity assumption we make the additional step of constructing confidence bands for $b$ that shrink at a rate adapting to the unknown smoothness.

Estimating the invariant density of a Markov process has been of interest for a long time. An early treatment is given by \citet{roussas1969}, who considered kernel estimators and showed consistency and asymptotic normality of the estimators under the strong Doeblin condition.
\citet{rosenblatt1970} analysed kernel estimators under the weaker condition $G_2$ on the Markov chain. More general $\delta$-sequences were used for the estimation by \citet{CastellanaLeadbetter1986}, who prove pointwise consistency and under strong mixing assumptions asymptotic normality.
\citet{yakowitz1989} shows asymptotic normality of kernel density estimators for the invariant density of Markov chains without using assumptions on the rates of mixing parameter sequences.
Adaptive estimation was considered by \citet{lacour2008}, who estimates the invariant density and the transition density of Markov chains by model selection and proves that the estimators attain the minimax convergence rate under $L^2$-loss.
For stationary processes \citet{schmisser2013} estimates the derivatives of the invariant density by model selection, derives the convergence rates of the estimators and pays special attention to the case of discretely observed diffusion processes. We see that asymptotic normality has been widely considered in the nonparametric estimation of invariant densities and thus implicitly also confidence intervals. However, we are not aware of any extensions of the pointwise results to uniform confidence bands for invariant densities, which are, for instance, necessary to construct goodness-of-fit tests of the Kolmogorov--Smirnov type.

The statistical properties of the diffusion model depend crucially on the observation
scheme. If the whole path $(X_t)_{0\le t\le T}$ is observed for some time
horizon $T>0$, we speak of \emph{continuous} observations. The case of discrete
observations $(X_{k\Delta})_{k=0,\dots,n-1}$ with observation distance 
$\Delta>0$
is distinguished into \emph{high-frequency} observations, i.e. $\Delta\downarrow
0$, and \emph{low-frequency} observations, where $\Delta>0$ is fixed. While in the first two settings path properties of the process can be used, statistical inference for low-frequency observations has to rely on the Markovian structure of the observations.
A review on parametric estimation in diffusion models is given by \citet{kutoyants2004} and \citet{aitsahalia2010}.
Nonparametric results are summarized in \citet{gobetEtAl2004}, where also 
estimators based on low-frequency observations are introduced and analysed. These low-frequency estimators rely on a spectral identification of diffusion coefficients which have been introduced by \citet{hansenScheinkman1995} and \citet{hansenEtAl1998}. On the same observation scheme \citet{kristensen2010} studies a pseudo-maximum likelihood approach in a semiparametric model. Nonparametric estimation based on random sampling times of the diffusion has been studied in \citet{chorowskiTrabs2015}. While we pursue a frequentist approach, the Bayesian approach is also very attractive. Based on low-frequency observations 
\citet{vanderMeulenvanZanten2013} have proved consistency of the Bayesian method and \citet{NicklSoehl2015} showed posterior contraction rates.

As usual, nonparametric estimators depend on some tuning parameters,
such as the bandwidth for classical kernel estimators. Choosing these parameters in a data-driven way, \citet{spokoiny2000} initiated adaptive drift estimation in
the diffusion model based on continuous observations. This was further
developed by \citet{Dalalyan2005} and \citet{Loecherbach2011}. Based on
high-frequency observations, adaptive estimation was studied by
\citet{hoffmann1999} as well as \citet{comteEtAl2007}. 
In the low-frequency case the question of adaptive estimation has been studied by \citet{chorowskiTrabs2015}. In this work we go one step further not only constructing a (rate optimal) adaptive estimator for the drift, but constructing
adaptive confidence bands.

Statistical applications require tests and confidence statements. 
\citet{NegriNishiyama2010} as well as \citet{masudaEtAl2011} have constructed 
goodness-of-fit tests for diffusions based on high-frequency observations. 
\citet{low1997} has shown that even in a simple density estimation problem no 
confidence bands exist which are honest and adaptive at the same time. 
Circumventing this negative result by a ``self-similarity'' condition, 
\citet{gineNickl2010} have constructed honest and adaptive confidence bands for density estimation. \citet{hoffmannNickl2011} have further studied necessary and sufficient conditions for 
the existence of adaptive confidence bands and the ``self-similarity'' condition has led to several recent papers on adaptive confidence bands, notably \citet{chernozhukovChetverikovKato2014} and 
\citet{SzabovanderVaartvanZanten2014}.
The present paper extends the theory of adaptive confidence bands beyond the classical nonparametric models of density estimation, white noise regression and the Gaussian sequence model which have been  treated in the above papers.

In order to derive confidence bands, we first have to establish a uniform
central limit theorem. The empirical measure of the observations
$X_0,\dots,X_{(n-1)\Delta}$ is the canonical estimator for the invariant measure
of a Markov chain or diffusion. Considering a wavelet projection estimator, we obtain a
smoothed version of the empirical measure, which is subsequently used to estimate the drift function in the case of diffusions. Thus a natural starting point is a functional central limit theorem 
for the invariant measure.
Since our observations are not independent, the standard empirical
process theory does not apply. Instead we have to use the Markov structure of
the chain $(X_{k\Delta})_k$. In the continuous time analogue the Donsker theorem
for diffusion processes has been studied by \citet{vanderVaartvanZanten2005}. In
the case of low-frequency observations, the estimation problem is ill-posed and 
we have nonparametric convergence rates under the uniform loss.
For the asymptotic behaviour of the estimation error in the uniform norm we 
would expect a Gumbel distribution as shown by 
\citet{gineNickl2010} in the density estimation case using extreme value 
theory.
Recent papers by \citet{castilloNickl2013,castilloNickl2014}
show that we can hope for parametric rates and an asymptotic normal
distribution if we consider instead
a weaker norm for the loss. More precisely, the estimation error can be measured
in a multi-scale space where the wavelet coefficients are down-weighted
appropriately. The resulting norm corresponds to a negative H\"older norm. 

Following this approach and relying on a concentration inequality by
\citet{adamczakBednorz2013}, our first result is a functional central limit 
theorem for rather general geometrically ergodic, Harris-recurrent Markov 
chains. This could also be of interest for the theory on Markov chain Monte Carlo 
(MCMC) methods considering that the central limit theorem measures the distance between a target integral
and its approximation,
\[
  \int_\R f(z)\mu(\d z)\quad\text{and}\quad \frac{1}{n}\sum_{k=0}^{n-1}f(Z_k),
\]
respectively, where $(Z_k)_k$ is a Markov chain with invariant measure $\mu$, cf. \citet{geyer1992}. Nevertheless, our focus is on the statistical point of view. The functional central limit theorem immediately yields non-adaptive confidence bands and as in 
\citet{castilloNickl2014} these have an $L^{\infty}$-diameter shrinking with
(almost) the optimal nonparametric rate.
This small deviation from the optimal rate corresponds to the
usual undersmoothing in the construction of nonparametric confidence sets.

Applying the results for general Markov chains to diffusion processes observed at low frequency, we obtain 
a functional central limit theorem for estimators of the drift function. Inspired 
by \citet{gineNickl2010}, in a last demanding step the smoothness of $b$ and the 
corresponding size of the confidence band is estimated to find adaptive 
confidence bands. The adaptive procedure relies on Lepski's method.
In order to make the construction of adaptive confidence bands feasible, we impose a self-similarity assumption on the drift 
function.

This work is organized as follows: In Section~\ref{sec:ucltMu} we study general
Markov chains and prove the functional central limit theorem and confidence bands
under appropriate conditions on the chain. These results are applied to
diffusion processes in Section~\ref{sec:AppDiff}. The adaptive confidence bands for
the drift estimator are constructed in Section~\ref{sec:AdaptiveConf}. Some
proofs are postponed to the last two sections.

\section{\label{sec:ucltMu}Confidence bands for the invariant probability
density of Markov processes}

\subsection{Preliminaries on Markov chains}

We start with recalling some facts from the theory of Markov chains. For all basic definition and results we refer to \citet{meynTweedie2009}. Let $Z=(Z_k)$, $k=0,1,\dots$, be a time-homogeneous Markov chain with state space $(\R,\mathcal{B}(\R))$. 
To fix the notation, let $P_x$ and $P_{\nu}$ denote the probability measure of the chain with initial conditions $Z_0=x\in\R$ and $Z_0\sim\nu$, respectively. The corresponding expectations will be denoted by $\E_x$ and $\E_\nu$, the Markov chain transition kernel by $P(x,A)$, $x\in\R$,
$A\in\mathcal{B}(\R)$. The transition operator is defined by $(Pf)(x)=\E_x[f(Z_1)]$. 

From the general theory of Markov chains we know that for a Harris-recurrent 
Markov chain~$Z$ the existence of a unique invariant probability measure~$\mu$ 
is equivalent to the \emph{drift condition}
\[
  PV(x)-V(x)\le -1+c\1_C(x)
\] 
for some petite set $C$, some $c<\infty$ and some non-negative function $V$, 
which is finite at some $x_0\in\R$. 
If $Z$ is additionally aperiodic, then this drift condition is already equivalent to $Z$ being \emph{ergodic} 
\[
  \|P^n(x,\cdot)-\pi\|_{TV}\to0,\quad \text{as}\quad n\to\infty,\quad\text{for all }x\in\R,
\] 
denoting the total variation norm of a measure by $\|\cdot\|_{TV}$.
If we impose a stronger drift condition, namely the geometric drift towards $C$, we obtain even \emph{geometric ergodicity}: 
For a $\psi$-irreducible and aperiodic Markov chain $Z$ satisfying
\begin{align}
 \label{con:drift}
(PV)(x)-V(x)\le-\lambda V(x)+c\1_C(x), \quad\text{ for all }x\in\R,
\end{align}
for a petite set $C$, some $\lambda>0,c<\infty$ and a function $V\colon\R\to[1,\infty)$, it holds for some $r>1$, $R<\infty$,
\[
  \sum_{n\ge0}r^n\|P^n(x,\cdot)-\mu\|_{TV}\le RV(x),\quad\text{for all
}x\in\R.
\]
Note that $\psi$-irreducibility together with the geometric drift condition 
\eqref{con:drift} implies already that~$Z$ is positive Harris with invariant 
probability measure~$\mu$. 

The geometric ergodicity yields the following central limit theorem, see \citet[Thm. II.4.1]{chen1999}. The weakest form of ergodicity so that the central limit theorem holds is ergodicity of degree~2 which is slightly weaker than the geometric ergodicity that we have assumed here.
\begin{prpstn}\label{prop:clt}
Let $(Z_k)_{k\ge0}$ be a geometrically ergodic Markov chain
with arbitrary initial condition and invariant probability measure $\mu$, 
then there exists for every bounded function $f=(f_1,\dots,f_d):\R\to\R^d$ a symmetric, positive semidefinite matrix
$\Sigma_f=(\Sigma_{f_i,f_j})_{i,j=1,\dots,d}$ such that
  \[ n^{-1/2}\Big(\sum_{k=0}^{n-1}f(Z_k)-n\E_\mu[f(Z_0)]\Big)\overset{d}{\longrightarrow} N(0,\Sigma_f),\quad\text{as }n\to\infty.
  \]
  For $i,j\in\{1,\dots,d\}$ the asymptotic covariances are given by
  \begin{align}
    \Sigma_{f_i,f_j} 
    &:=\lim_{n\to\infty}n^{-1} \Cov_{\mu}\Big(\sum_{k=0}^{n-1}f_i(Z_{k}),\sum_{k=0}^{n-1}f_j(Z_{k})\Big)\label{eq:Variance}\\
    &\phantom{:}=\E_{\mu}\big[(f_i(Z_{0})-\E_{\mu}[f_i])(f_j(Z_{0})-\E_{\mu}[f_j])\big] +\sum_{k=1}^{\infty}\E_{\mu}\big[(f_i(Z_{0})-\E_{\mu}[f_i])(f_j(Z_{k})-\E_{\mu}[f_j])\big]\nonumber\\
    &\qquad+\sum_{k=1}^{\infty}\E_{\mu}\big[(f_i(Z_{k})-\E_{\mu}[f_i])(f_j(Z_{0})-\E_{\mu}[f_j])\big].\nonumber 
  \end{align}
\end{prpstn}

In order to lift this ``pointwise'' result to a functional central limit
theorem, we will in addition need a concentration inequality for a preciser
control on how the sum $n^{-1}\sum_{k=0}^{n-1}f(Z_{k})$ deviates from the
integral $\int f(z)\mu(\d z)$ for finite
sample sizes.
To this end, we strengthen the aperiodicity assumption to
\emph{strong aperiodicity} \citep[see][Prop. 5.4.5]{meynTweedie2009}, that is
there exists a set $C\in\mathcal{B}(\R)$, a probability measure $\nu$ with
$\nu(C)>0$ and a
constant $\delta>0$ such that
\begin{align} \label{con:smallset}
P(x,B)\ge\delta\nu(B), \quad \text{ for all }x\in C, B\in \mathcal{B}(\R).
\end{align}
Any set $C$ satisfying this condition is called \emph{small set}. Recall that any small set is a petite set.

\begin{prpstn}[Theorem~9 by \citet{adamczakBednorz2013}]\label{prop:Concentration}
  Let $Z=(Z_k)_{k\ge0}$ be a Harris recurrent, strongly aperiodic Markov chain on $(\R,\mathcal{B}(\R))$ with unique invariant measure $\mu$. For some set $C\in\mathcal B(\R)$ with $\mu(C)>0$ let $Z$ satisfy the drift condition~\eqref{con:drift} and the small set condition~\eqref{con:smallset}.

  Let $f\in L^{2}(\mu)$ be bounded. For any $0<\tau<1$ there are constants 
$K,c_{2}$ depending only on $\delta$, $V$, $\lambda$, $c$ and $\tau$ and a constant 
$c_{1}$ depending additionally on the initial value $x\in\R$ such that for any 
$t>0$
  \begin{align*}
    & P_{x}\Big(\Big|\sum_{k=0}^{n-1}f(Z_{k})-n\E_{\mu}[f(Z_0)]\Big|>t\Big)\\
    & \qquad\le K\exp\left(-c_{1}\Big(\frac{t}{\|f\|_{\infty}}\Big)^{\tau}\right)+K\exp\bigg(-\frac{
  c_{2}t^{2}}{n\Sigma_{f}+t\max(\|f\|_{\infty}(\log
  n)^{1/\tau},\Sigma_{f}^{1/2})}\bigg),
  \end{align*}
  where $\Sigma_f$ is given by~\eqref{eq:Variance} with $d=1$.
\end{prpstn}

As a last ingredient we need to bound the asymptotic variance $\Sigma_{f}$
in Propositions~\ref{prop:clt} and \ref{prop:Concentration} in terms of $\|\bar f\|_{L^{2}(\mu)}^{2}$ for the centred function $\bar f:=f-\int f \d\mu$. The geometric ergodicity only yields a bound $\mathcal{O}(\|\bar f\|_{\infty}^{2})$. Therefore, we require that the transition operator is a contraction in the sense that there exists some $\rho\in(0,1)$ satisfying
\begin{align}\label{con:contraction}
 \|Pg\|_{L^{2}(\mu)}\le\rho\|g\|_{L^{2}(\mu)}\quad\text{ for
all }g\in L^2(\mu)\text{ with }\int g\,\d\mu=0.
\end{align}
This property is also known as $\rho$-mixing. It corresponds to a Poincar\'e inequality \citep[cf.][Thm. 1.3]{bakryEtAl2008} and its relation to drift conditions is analysed by \citet{bakryEtAl2008}. If \eqref{con:contraction} is fulfilled, the Cauchy--Schwarz inequality yields
\begin{equation}
\Sigma_{f}=\Sigma_{\bar f}\le\|\bar f\|_{L^{2}(\mu)}^{2} +2\sum_{k=1}^{\infty}\|\bar f\|_{L^{2}(\mu)}\|P^{k}\bar f\|_{L^{2}(\mu)}\le\big(1+2\sum_{k=1}^{\infty}\rho^{k}\big)\|\bar f\|_{L^{2}(\mu)}^{2}= \tfrac{1+\rho}{1-\rho}\|\bar f\|_{L^{2}(\mu)}^{2}.\label{eq:estSigma}
\end{equation}

\subsection{A functional central limit theorem}\label{sec:UCLTMarkov}
The basic idea is to prove a functional central limit theorem for the invariant
probability measure~$\mu$ by choosing an orthonormal basis, applying the 
pointwise
central limit theorem to the basis functions (Proposition~\ref{prop:clt}) and
extending this result to finite linear combinations with the help of the
concentration inequality (Proposition~\ref{prop:Concentration}). Provided~$\mu$ 
has some
regularity, the approximation error due to considering only a finite basis
expansion of $\mu$ will be negligible. Noting that it is straightforward to 
extend the results to any compact subset of $\R$, we focus on a central limit 
theorem on a bounded interval $[a,b]$ with $-\infty<a<b<\infty$.

Let $(\phi_{j_0,l},\psi_{j,k}:j\ge j_0,l,k\in\Z)$, for some $j_0\ge0$, a scaling function $\phi$ and a wavelet function $\psi$, be a regular compactly supported $L^{2}$-orthonormal wavelet basis of $L^{2}(\R)$. For the sake of clarity we throughout use Daubechies' wavelets of order $N\in\N$, but any other compactly supported regular wavelet basis can be applied as well. As a standing assumption we suppose that $N$ is chosen large enough such that the H\"older regularity of $\phi$ and $\psi$ is larger than the regularity required for the invariant measure. The approximation spaces for resolution levels $J>j_0$ are defined as
\[
V_{J}:=\overline{\Span}\{\phi_{j_0,l},\psi_{j,k}:j=j_0,\dots,J,\;l,k\in \Z\},
\]
The projection onto $V_{J}$ is denoted by $\pi_{J}$. Since $j_0$ is fixed and to simplify the notation, we write $\psi_{-1,l}:=\phi_{j_0,l}$.

Using the first $n\in\N$ steps $Z_{0},Z_{1},\dots,Z_{n-1}$ of a realisation of the chain,
we define the empirical measure 
\[
\mu_{n}:=\frac{1}{n}\sum_{k=0}^{n-1}\delta_{Z_{k}},
\]
where $\delta_{x}$ denotes the Dirac measure at the point $x\in\R$. The 
canonical projection wavelet estimator of $\mu$ is given by 
\begin{equation}
  \hat{\mu}_{J}:=\pi_{J}(\mu_{n})
  =\sum_{l\in\Z}\hat \mu_{-1,l}\psi_{-1,l}+\sum_{j=j_0 }^J\sum_{k\in\Z}\hat{\mu}_{j,k}\psi_{j,k},\quad\quad\hat{\mu}_{
  j,k}:=\langle\psi_{j,k},\mu_{n}\rangle:=\int\psi_{j,k}\,\d\mu_{n}.
  \label{eq:hatMu}
\end{equation}
For any $\psi_{j,k}$ Proposition~\ref{prop:clt} yields that
$\sqrt{n}(\mu_{n}-\mu)(\psi_{j,k})$
converges in distribution for $n\to\infty$ to a Gaussian random variable
\begin{gather}
\mathbb{G_{\mu}}(j,k)\sim\mathcal{N}(0,\Sigma_{\psi_{j,k}})
\quad\text{with covariances }\quad \E_{\mu}[\mathbb{G}_{\mu}({j,k})\mathbb{G}
_{\mu}({l,m})]= \Sigma_{\psi_{j,k},\psi_{l,m}}.\label{eq:limit}
\end{gather}
Using the techniques from \citet{castilloNickl2014}, this pointwise
convergence of $\mu_{n}$ can be extended to a uniform central limit
theorem on $[a,b]$ for the projection estimator $\hat{\mu}_{J}$ 
in the \emph{multi-scale sequence spaces} which are defined as follows: 
Noting that the Daubechies wavelets fulfil $\supp\phi\subset[0,2N-1]$ and 
$\supp\psi\subset[-N+1,N]$, cf. \citet[Chap. 7]{hardleEtAl1998}, 
the sets $L:=K_{-1}:=\{k\in\Z:2^{j_0}a-2N+1\le k\le 2^{j_0}b\}$ and $K_j:=\{k\in\Z:2^ja-N\le k\le 2^jb+N-1\}$ contain all indices of $\phi_{j_0,\cdot}$ and $\psi_{j,\cdot}$, respectively, whose support intersects with the interval $[a,b]$. For a monotonously increasing weighting sequence $w=(w_{j})_{j=-1,j_0,j_0+1,j_0+2,\dots}$ with $w_{j}\ge1$ and $w_{-1}:=1$ we define the multi-scale sequence spaces as
\[
  \mathcal{M}:=\mathcal{M}(w):=\Big\{x=(x_{jk}):\|x\|_{\mathcal{M}(w)}:=\sup_{j\in\{-1,j_0,j_0+1,\dots\}}\max_{k\in K_j}\frac{|x_{jk}|}{w_{j}}<\infty\Big\},
\] 
Since the Banach space $\mathcal{M}(w)$ is non-separable, we define
the separable, closed subspace 
\[
\mathcal{M}_{0}:=\mathcal{M}_{0}(w):=\Big\{
x=(x_{jk}):\lim_{j\to\infty}\max_{k\in K_j}\frac{|x_{jk}|}{w_{j}}=0\Big\}.
\]

Let us assume that $\mu$ is absolutely continuous with respect to the Lebesgue
measure and denote the density likewise by $\mu$.
If the density is bounded on $D=[a-2^{-j_0}(2N-1),b+2^{-j_0}(2N-1)],$ the orthonormality and the support of $(\psi_{j,k})$ and
(\ref{eq:estSigma}) yield 
$\Sigma_{\psi_{j,k}}=\mathcal{O}(\|\mu\|_{\infty})$. Standard
estimates of the supremum of normal random variables yield that the maximum over 
the $2^j$ variables $\mathbb G_\mu({j,\cdot})$ of a resolution level $j$ are of 
the order $\max_k|\mathbb G_\mu({j,k})|=\mathcal O_P(\sqrt j)$, see 
\eqref{eq:maxGauss} below. Since the cardinality of $K_j$ is of the order 
$2^j$, a weighting $w_j=\sqrt j$ seems to be
appropriate and indeed we conclude as \citet[Prop. 3]{castilloNickl2014}:
\begin{lmm}
\label{lem:PropLimit}
Let $\mu$ admit a Lebesgue density which is bounded on $D$.
Then $\mathbb{G_{\mu}}$ from (\ref{eq:limit}) satisfies
$\E[\|\mathbb{G}_{\mu}\|_{\mathcal{M}(w)}]<\infty$ for the weights $w_{j}=\sqrt{j}$. Moreover, $\mathcal{L}(\mathbb{G}_{\mu})$
is a tight Gaussian Borel probability measure in $\mathcal{M}_{0}(w)$
if $\sqrt{j}/w_{j}\to0$. 
\end{lmm}

Let us now summarise the assumptions on the Markov chain, which are needed to prove the functional central limit theorem and for the construction of confidence bands. For any regularity $s>0$, denoting the integer part of $s$ by $[s]$, the \emph{H\"older space} on a domain $D$ is defined by
\[
  C^s(D):=\Big\{f\colon D\to\R\Big|\|f\|_{C^s}:=\sum_{k=0}^{[s]}\|f^{(k)}\|_\infty+\sup_{x\neq y}\frac{|f^{[s]}(x)-f^{[s]}(y)|}{|x-y|^{s-[s]}}\Big\}.
\]
\begin{assumption}\label{ass:Markovchain}
  Let $(Z_k)_{k\ge0}$ be a Harris recurrent, strongly aperiodic Markov chain on 
$(\R,\mathcal B(\R))$ with initial condition $Z_0=x$. Let the invariant 
probability measure have a density $\mu$ in $C^{s}(D)$ for some $s>0$ and some 
sufficiently large set $D\subset\R$ containing $[a,b]$. Let the drift
  condition~\eqref{con:drift} and small set condition~\eqref{con:smallset} be satisfied for some $C\in\mathcal B(\R)$ with $\mu(C)>0$. Further suppose that the transition operator is an $L^2(\mu)$-contraction fulfilling \eqref{con:contraction} with $\rho\in(0,1)$.
\end{assumption}
\begin{rmrk}
  As we have discussed above it suffices to verify that the chain $(Z_k)_{k\ge0}$ is $\psi$-irreducible and satisfies \eqref{con:drift} and \eqref{con:smallset} in order to conclude that the $(Z_k)_{k\ge0}$ is Harris recurrent, strongly aperiodic and has a unique invariant probability measure.
\end{rmrk}

Now we can show the functional central limit theorem for $\hat{\mu}_{J}$ in the
space $\mathcal{M}_{0}(w)$. Note that the natural nonparametric choice $J_{n}$ given by $2^{J_{n}}\sim n^{1/(2s+1)}$ satisfies the conditions of the following theorem. Recall that weak convergence of laws $\mathcal{L}(X)$ of random variables $X$
on a metric space $(S,d)$ can be metrised by the \emph{bounded-Lipschitz
metric}
\begin{align*}
\beta_{S}(\mu,\nu) & :=\sup_{F:\|F\|_{BL}\le1}\Big|\int_{S}F(x)(\mu(\d x)-\nu(\d
x))\Big|\quad\text{with}\\
\|F\|_{BL} & :=\sup_{x\in S}|F(x)|+\sup_{x,y\in S:x\neq
y}\frac{|F(x)-F(y)|}{d(x,y)}.
\end{align*}

\begin{thrm}\label{thrm:ucltMu}
  Grant Assumption~\ref{ass:Markovchain} and let $w=(w_{j})$ be increasing and
  satisfy $\sqrt{j}/w_{j}\to0$ as $j\to\infty$. Let $J_{n}\in\N$ fulfil, for some $\tau\in(0,1)$,
  \[
  \sqrt{n}2^{-J_{n}(2s+1)/2}w_{J_{n}}^{-1}=o(1),\qquad(\log
n)^{2/\tau}n^{-1}2^{J_{n}}J_{n}=\mathcal{O}(1).
  \]
  Then $\hat{\mu}_{J_{n}}$ from (\ref{eq:hatMu}) satisfies, for $n\to\infty$,
  \[ 
\sqrt{n}(\hat{\mu}_{J_{n}}-\mu)\overset{d}{\longrightarrow}\mathbb{G}_{\mu}
\quad\text{in }\mathcal{M}_{0}(w).
  \]
\end{thrm}
\begin{proof}
We follow the strategy of \cite[Thm. 1]{castilloNickl2014}. First
we deal with the bias term. By the $s$-H\"older regularity of
$\mu$ we have \cite[Definition (5.90) and Proposition 5.3.13]{gineNickl2015}
\[
\sup_{j,k}2^{j(2s+1)/2}|\langle\psi_{j,k},\mu\rangle|<\infty
\]
and thus by the assumption on $J_{n}$
\[
\|\mu-\pi_{J_{n}}(\mu)\|_{\mathcal{M}}=\sup_{j>J_{n}}\max_{k\in K_j}w_{j}^{-1}
|\langle\psi_{jk},\mu\rangle|\lesssim\sup_{j>J_{n}}w_{j}^{-1}2^{-j(2s+1)/2}=o(n^
{-1/2}).
\]
Defining $\nu_{n}:=\sqrt{n}(\hat{\mu}_{J_{n}}-\pi_{J_{n}}(\mu))$, we
decompose the stochastic error, for $J<J_{n}$ to be specified later,
\begin{align}
\beta_{\mathcal{M}_{0}}(\mathcal{L}(\nu_{n}),\mathcal{L}(\mathbb{G}_{\mu})) &
\le\beta_{\mathcal{M}_{0}}(\mathcal{L}(\nu_{n}),\mathcal{L}(\nu_{n})\circ\pi_{J}
^{-1})+\beta_{\mathcal{M}_{0}}(\mathcal{L}(\nu_{n})\circ\pi_{J}^{-1},\mathcal{L}
(\mathbb{G}_{\mu})\circ\pi_{J}^{-1})\nonumber \\
 &
\qquad+\beta_{\mathcal{M}_{0}}(\mathcal{L}(\mathbb{G}_{\mu})\circ\pi_{J}^{-1},
\mathcal{L}(\mathbb{G}_{\mu})).\label{eq:ucltDecomp}
\end{align}
In the sequel we will separately show that all three terms converge
to zero. Let $\epsilon>0$. By definition of the $\beta_{\mathcal{M}_{0}}$-norm
we estimate the first term by 
\begin{align}
  \beta_{\mathcal{M}_{0}}(\mathcal{L}(\nu_{n}),\mathcal{L}(\nu_{n})\circ\pi_{J}^{-1})
  &=\sup_{F:\|F\|_{BL}\le1}\big|\E\big[F(\nu_n)-F(\pi_J(\nu_n))\big]\big|\nonumber\\
  &\le \E\big[\|\sqrt{n}(\pi_{J_{n}}-\pi_{J})(\mu_{n}-\mu)\|_{\mathcal{M}}\big]\nonumber \\
  &\le  \max_{J<j\le J_{n}}(w_{j}^{-1}j^{1/2})\E\Big[\max_{J<j\le J_{n}}\max_{k\in K_j}j^{-1/2}|\langle\sqrt{n}(\mu_{n}-\mu),\psi_{j,k}\rangle|\Big].\label{eq:ucltFirst}
\end{align}
By the assumptions on $w$ and due to the factor in front of the expectation,
the above display can be bounded by $\epsilon/3$ if $J$ is chosen
large enough and provided that the expectation can be bounded by a
constant independent of $J$ and $n$. To apply the concentration
inequality in Proposition~\ref{prop:Concentration}, note that
$\Sigma_{\psi_{j,k}}=\mathcal{O}(\|\mu\|_{\infty})$
by (\ref{eq:estSigma}) and
$\sqrt{j}\|\psi_{j,k}\|_{\infty}=\sqrt{j}2^{j/2}=\mathcal{O}(\sqrt{n}(\log
n)^{-1/\tau})$
for $j\le J_{n}$. Hence, for any $M>0$ large enough we obtain for
constants $c_{i}>0,i=1,2,\dots,$
\begin{align}
 & \quad\E\Big[\max_{J<j\le
J_{n}}\max_{k\in K_j}j^{-1/2}|\langle\sqrt{n}(\mu_{n}-\mu),\psi_{j,k}\rangle|\Big]
\nonumber \\
&\le  M+\int_{M}^{\infty}P\Big(\max_{J<j\le
J_{n}}\max_{k\in K_j}j^{-1/2}|\langle\sqrt{n}(\mu_{n}-\mu),\psi_{j,k}\rangle|>u\Big)\d
u\nonumber \displaybreak[0]\\
&\le  M+\sum_{J<j\le
J_{n},k\in K_j}\int_{M}^{\infty}P\Big(|\langle\sqrt{n}(\mu_{n}-\mu),\psi_{j,k}
\rangle|>\sqrt{j}u\Big)\d u\nonumber \\
&\lesssim  M+\sum_{J<j\le J_{n}}2^{j}\int_{M}^{\infty}\Big(\exp\big(-c_{1}(\log
n)j^{\tau}u^{\tau}\big)+\exp\big(-c_{2}ju^{2}/(1+u)\big)\Big)\d
u\label{eq:ApplyBernstein}\\
&\lesssim  M+\sum_{J<j\le J_{n}}2^{j}\Big(\frac{e^{-c_{3}(jM)^{\tau}\log
n}}{j^{\tau}\log n}+\frac{e^{-c_{4}jM}}{j}\Big)\nonumber \\
&\lesssim  M+e^{-c_{5}JM^{\tau}}\lesssim M+1,\nonumber 
\end{align}
where we have used in the next to last estimate that $J_n\lesssim \log n$ and thus
$j^\tau\log n\gtrsim j$ for all $j\le J_n$.

To bound the second term in (\ref{eq:ucltDecomp}), we use
Proposition~\ref{prop:clt}
and the Cram\'{e}r--Wold device to see that it is smaller than $\epsilon/3$
for fixed $J$ and $n$ sufficiently large. It remains to consider
the third term in (\ref{eq:ucltDecomp}) which can be estimated similarly
to (\ref{eq:ucltFirst}), using that
$\E[\sup_{j}\max_{k}j^{-1/2}|\mathbb{G}_{\mu}({j,k})|]<\infty$
by Lemma~\ref{lem:PropLimit}.
\end{proof}

\subsection{The construction of confidence bands}\label{sec:ConstConfBand}

Using the multi-scale central limit theorem, we now construct confidence
bands for the density of the invariant probability measure. For some confidence
level $\alpha\in(0,1)$ the natural idea is to take 
\[
{\mathcal C}_{n}(\zeta_{\alpha}):=\Big\{ f:\|f-\hat{\mu}_{J_{n}}\|_{\mathcal M}<\frac{\zeta_{\alpha}}{\sqrt{n}}\Big\}
=\Big\{ f:\sup_{j,k}w_{j}^{-1}|\langle f-\hat{\mu}_{J_{n}},\psi_{j,k}\rangle|<\frac{\zeta_{\alpha}}{\sqrt{n}}\Big\},
\]
where $\zeta_{\alpha}$ is chosen such that
$P(\|\mathbb{G}_{\mu}\|_{\mathcal{M}}<\zeta_{\alpha})\ge1-\alpha$.
For this set the asymptotic coverage follows immediately from
Theorem~\ref{thrm:ucltMu}.
However, ${\mathcal C}_{n}(\zeta_{\alpha})$ is too large in terms of the $L^{\infty}([a,b])$-diameter
\[
|{\mathcal C}_{n}(\zeta_{\alpha})|_{\infty}:=\sup\big\{\sup_{x\in[a,b]}|f(x)-g(x)|:f,g\in {\mathcal C}_{n}(\zeta_{\alpha})\big\}.
\]
To obtain the (nearly) optimal $L^{\infty}$-diameter, we need to
control the large resolution levels. As suggested by \citet{castilloNickl2014},
we use a-priori knowledge of the regularity $s$ to define 
\begin{equation}
\bar{\mathcal C}_{n}:=\bar{\mathcal C}_{n}(\zeta_{\alpha},s,u_{n})
:={\mathcal C}_{n}(\zeta_{\alpha})\cap\big\{ f:\|f\|_{C^{s}}\le u_{n}\big\}\label{eq:CS}
\end{equation}
for a sequence $u_{n}\to\infty$.
\begin{prpstn}\label{prop:CSnonadaptive}
  Grant Assumption~\ref{ass:Markovchain} with
  $s>0$ and let $w=(w_{j})$ be increasing and satisfy $\sqrt{j}/w_{j}\to0$. For $\alpha\in(0,1)$ let $\zeta_{\alpha}>0$ be such that $P(\|\mathbb{G}_{\mu}\|_{\mathcal{M}}\ge\zeta_{\alpha})\le\alpha$ and choose $J_{n}:=J_{n}(s)$ such that 
  \[
  2^{J_{n}}=\Big(\frac{n}{\log
n}\Big)^{1/(2s+1)}.
  \]
  Then the confidence set $\bar{\mathcal C}_{n}=\bar{\mathcal C}_{n}(\zeta_{\alpha},s,u_{n})$
  from (\ref{eq:CS}) with $u_{n}:=w_{J_{n}}/\sqrt{J_{n}}$ satisfies
  \[ 
    \liminf_{n\to\infty}P(\mu\in\bar{\mathcal C}_{n})\ge1-\alpha\quad\text{and}\quad|\bar{\mathcal C}_{n}|_{\infty}
    =\mathcal {O}_{P}\Big(\Big(\frac{n}{\log n}\Big)^{-s/(2s+1)}u_{n}\Big).
  \]
\end{prpstn}
\begin{proof}
Let us first verify $\liminf_{n\to\infty}P(\mu\in\bar{\mathcal C}_{n})\ge1-\alpha.$ Since $\mu\in
C^{s}(u_{n})$
for large enough $n$, Theorem \ref{thrm:ucltMu} yields 
\begin{align*}
\liminf_{n\to\infty}P(\mu &
\in\bar{\mathcal C}_{n})
=\liminf_{n\to\infty}P(\sqrt{n}\|\hat{\mu}_{J_{n}}-\mu\|_{
\mathcal
{M}}<\zeta_{\alpha})\ge P(\|\mathbb{G}_{\mu}\|_{\mathcal{M}}<\zeta_{\alpha}
)\ge1-\alpha.
\end{align*}
To bound the diameter let $f,g\in\bar{\mathcal C}_{n}$. Using
$\|f-\hat{\mu}_{J_{n}}\|_{\mathcal{M}}=\mathcal{O}_{P}(n^{-1/2})$
and $f-g\in C^{s}(2u_{n})$, we obtain
\begin{align}
\|f-g\|_{L^\infty([a,b])}&\lesssim  \sum_{j\le J_{n}}2^{j/2}\max_{k\in K_j}|\langle
f-g,\psi_{j,k}\rangle|+\sum_{j>J_{n}}2^{j/2}\max_{k\in K_j}|\langle
f-g,\psi_{j,k}\rangle|\nonumber \\
&\le  \sum_{j\le J_{n}}2^{j/2}\big(\max_{k\in K_j}|\langle
f-\hat{\mu}_{J_{n}},\psi_{j,k}\rangle|+\max_{k\in K_j}|\langle
g-\hat{\mu}_{J_{n}},\psi_{j,k}\rangle|\big)\nonumber \\
 & \qquad+\sum_{j>J_{n}}2^{-js}2^{j(s+1/2)}\max_{k\in K_j}|\langle
f-g,\psi_{j,k}\rangle|\nonumber \\
&\le 
\big(\|f-\hat{\mu}_{J_{n}}\|_{\mathcal{M}}+\|g-\hat{\mu}_{J_{n}}\|_{\mathcal{M}}
\big)\sum_{j\le J_{n}}2^{j/2}w_{j}+\|f-g\|_{C^{s}}\sum_{j>J_{n}}2^{-js}\nonumber
\\
&=
\mathcal{O}_{P}\Big(n^{-1/2}2^{J_{n}/2}w_{J_{n}}\Big)+\mathcal{O}_{P}\Big(2^{-J_
{n}s}u_{n}\Big)\nonumber \\
&= 
\mathcal{O}_{P}\Big(n^{-1/2}2^{J_{n}/2}J_{n}^{1/2}u_{n}\Big)+\mathcal{O}_{P}
\Big(2^{-J_{n}s}u_{n}\Big).\label{eq:diameter}
\end{align}
Plugging in the choice of $J_{n}$, we finally have
$n^{-1/2}\sqrt{2^{J_{n}}J_{n}}\lesssim(n/\log n)^{-s/(2s+1)}=2^{-J_{n}s}$.
\end{proof}

A multi-scale confidence band as in \eqref{eq:CS} allows for the construction of a classical $L^\infty$-band on $[a,b]$ around $\hat\mu_{J_n}$ as follows: Let us denote the almost optimal diameter by $\rho_n:=\big(\frac{n}{\log n}\big)^{-s/(2s+1)}u_{n}$. As we can deduce from \eqref{eq:diameter}, there is a constant $D>0$ such that we have $\|f-\hat\mu_{J_n}\|_{L^\infty([a,b])}\le D\rho_n$ for any $f\in \bar{\mathcal C_n}$. Hence, the band
\[
  \tilde{\mathcal C}_n:=\big\{f\colon [a,b]\to \R\big|\|f-\hat\mu_{J_n}\|_{L^\infty([a,b])}\le D\rho_n \big\}
\]
contains $\bar{\mathcal C}_n$ which only improves the coverage. Consequently, $\tilde{\mathcal C}_n$ is an $L^\infty$-confidence band with level $\alpha$ which shrinks with almost optimal rate $\rho_n$. In addition, a multi-scale confidence band as in \eqref{eq:CS} allows for simultaneous confidence intervals in all wavelet coefficients. This is especially useful for goodness-of-fit tests where the optimal $L^\infty$-diameter of $\bar{\mathcal C}_n$ is a measure of the power of the test.

In order to apply the confidence band \eqref{eq:CS} we need the regularity $s$ of the invariant density and a critical value $\zeta_{\alpha}$ such that
$P(\|\mathbb{G}_{\mu}\|_{\mathcal{M}}<\zeta_{\alpha})\ge1-\alpha$ for
$\alpha\in(0,1)$. 
Adaptive confidence bands will be presented later in the context of diffusions. 
So let us suppose for a moment that the regularity $s$ is known. Then the
problem reduces to the construction of the critical value to which the
remainder of this section is devoted.

A first observation is that
if several independent copies of the
diffusion are observed then one could calculate for each copy an estimator
$\hat\mu_{J_n}$ and obtain estimators for the values $\zeta_\alpha$ from
the distribution of the estimators $\hat\mu_{J_n}$ around their joint mean.
Since the assumption of many independent copies is not realistic we will not pursue
this further. Instead of the consistent estimation of the lowest possible $\zeta_\alpha$ we restrict ourselves to estimating an upper bound, which yields possibly more conservative confidence sets. By the concentration of Gaussian measures we know
for any $\kappa>0$ that
\[
  P\big(\|\mathbb{G}_{\mu}\|_{\mathcal{M}}\ge
\E[\|\mathbb{G}_{\mu}\|_{\mathcal{M}}]+\kappa\big)\le e^{-\kappa^2/(2\Sigma)},
\]
where $\Sigma:=\sup_{j,k}\E[|\mathbb G_\mu(j,k)|^2]=\sup_{j,k}\Sigma_{\psi_{j,k}}$, see for example \citet[Thm. 7.1]{ledoux2001}. Hence, an upper bound for $\zeta_\alpha$ is given by
\[
  \sqrt{2\Sigma\log \alpha^{-1}}+\E[\|\mathbb{G}_{\mu}\|_{\mathcal{M}}].
\]
The expected value $\E[\|\mathbb{G}_{\mu}\|_{\mathcal{M}}]$ can be bounded as in Proposition~2 by \citet{castilloNickl2014}, depending on $\Sigma$ again. We obtain the following upper bound for $\zeta_\alpha$:
\begin{lmm}
  Let $j_0\ge1$ and $w=(w_j)$ satisfy $w_{-1}=\sqrt{j_0}$ and $\inf_jw_j/\sqrt j\ge1,j\ge j_0,$ and define $\Sigma:=\sup_{j,k}\Sigma_{\psi_{j,k}}$. Then
  $P(\|\mathbb{G}_{\mu}\|_{\mathcal{M}}\ge\bar\zeta_{\alpha})\le\alpha$ holds for \begin{equation}
    \bar\zeta_\alpha(\Sigma):=\big(\sqrt{2\log 
\alpha^{-1}}+2C+\tfrac{32}{3C}2^{-2j_0}\big)\sqrt\Sigma\label{eq:barZeta}
  \end{equation}
  with $C:=(\sup_{j\ge j_0}(4\log|K_j|+2\log 2)/j)^{1/2}$.
\end{lmm}
\begin{proof}
  The cardinality of $K_j$ is denoted by $|K_j|$. Recall that a standard normal 
random variable~$Z$ satisfies 
  $$\E[e^{Z^2/4}]=\sqrt 2\quad\text{and}\quad P(Z>\kappa)\le\frac{1}{\kappa\sqrt{2\pi}}e^{-\kappa^2/2},\kappa>0.$$
  For each $j\ge j_0$ and $\kappa=2\sup_k\Sigma^{1/2}_{\psi_{j,k}}$ Jensen's inequality thus yields
  \begin{align}
    \E[\max_k|\mathbb G_\mu(j,k)|]
    &\le \kappa\Big(\log\E\big[e^{\max_k|\mathbb G_\mu(j,k)|^2/\kappa^2}\big]\Big)^{1/2}\notag\\
    &\le 2\sup_k\Sigma^{1/2}_{\psi_{j,k}}\big(\log|K_j|+\tfrac12\log2 \big)^{1/2}
    \le C\sqrt{\Sigma j},\label{eq:maxGauss}
  \end{align}
  for the constant $C:=(\sup_j(4\log|K_j|+2\log 2)/j)^{1/2}$. Theorem 7.1 in \cite{ledoux2001} yields for all~$t,T$ such that $t\ge \Sigma^{1/2} CT$ and $T>1$
  \begin{align*}
    P\big(\|\mathbb G_\mu\|_{\mathcal M}>t\big)
    &\le \sum_jP\big(|\max_k\mathbb G_\mu(j,k)-\E[\max_k\mathbb 
G_\mu(j,k)]|>tw_j-\E[\max_k|\mathbb G_\mu(j,k)|]\big)\\
    &\le \sum_jP\Big(|\max_k\mathbb G_\mu(j,k)-\E[\max_k\mathbb 
G_\mu(j,k)]|>(t-C)\sqrt {\Sigma j}\Big)\\
    &\le2\sum_je^{-jt^2(T-1)^2/(2\Sigma T^2)}. 
  \end{align*}
  Recall that $j=j_0,j_0,j_0+1,j_0+2,\dots$ in the above sum. Using Fubini's
theorem and the Gaussian tail bound, we conclude
  \begin{align*}
    \E[\|\mathbb G_\mu\|_{\mathcal M}]
    &\le \Sigma^{1/2} CT+\int_{\Sigma^{1/2} CT}^\infty P(\|\mathbb G_\mu\|_{\mathcal M}> t)\d t
    \le \Sigma^{1/2} CT+2\Sigma^{1/2}\sum_j\int_{CT}^\infty e^{-jt^2(1-1/T)^2/2}\d t\\
    &\le \Sigma^{1/2} CT+\frac{2\Sigma^{1/2} T}{C(T-1)^2}\sum_j e^{-(2\log 2)j(T-1)^2}
    \le \Sigma^{1/2} CT+\frac{4\Sigma^{1/2} T2^{-2(T-1)^2j_0}}{C(T-1)^2(1-2^{-2(T-1)^2})}.
  \end{align*}
  Choosing the $T=2$, we obtain $\E[\|\mathbb G_\mu\|_{\mathcal 
M}]\le(2C+\tfrac{32}{3C}2^{-2j_0})\Sigma^{1/2}$.
\end{proof}

From the above lemma we see that $\Sigma$ is the key quantity for the
construction of the critical values $\zeta_\alpha$. 
A natural estimator for $\Sigma$ is $\hat\Sigma_n:=(\max_{j\le
J_n,k}\hat\Sigma_{\psi_{j,k}})$, where $\hat\Sigma_{\psi_{j,k}}$ are 
estimators of~$\Sigma_{\psi_{j,k}}$ based on~$n$ observations. Since $J_n$ tends to infinity, 
the
maximum over all $j\le J_n$ converges to the supremum over all $j$ so that we
are asymptotically estimating the right quantity. 
For the estimators $\hat\Sigma_{\psi_{j,k}}$ we propose the initial monotone
sequence estimators based on autocovariations by \citet{geyer1992}, which are
consistent over-estimates, and this yields
almost surely 
\[\liminf_{n\to\infty}\hat\Sigma_n\ge\Sigma,\]
which suffices for our purposes.

The estimation of $\Sigma_{\psi_{j,k}}$ amounts to the estimation 
of the asymptotic variance $\Sigma_{f}$ in
(\ref{eq:Variance}) for a known function $f$ and this problem is studied in the
MCMC-literature. In addition to the sequence estimators, \citet{geyer1992}
discusses two other constructions together with their advantages and
disadvantages.
\citet{robert1995} constructs another estimator applying renewal theory, which
is however difficult to calculate. A more recent estimator using i.i.d.
copies of the process $X$ is given by \citet{chauveauDieboldt2003}.

As an alternative to the above estimation of $\Sigma$ in~\eqref{eq:barZeta} an upper bound could be estimated as follows:
Using~\eqref{eq:estSigma} we can bound $\Sigma$ from
above,
\[
  \Sigma\le \sup_{j,k}\tfrac{1+\rho}{1-\rho}\|\bar\psi_{j,k}\|_{L^2(\mu)}
  \le \sup_{j,k}\tfrac{1+\rho}{1-\rho}\|\psi_{j,k}\|_{L^2}\|\mu\|_\infty=\tfrac{1+\rho}{1-\rho}\|\mu\|_\infty,
\]
where we can plug in estimators for $\|\mu\|_\infty$ and $\rho$. 
Considering a wavelet $\psi_{j,k}$ localised around the maximum of $\mu$ we see
that the second inequality should provide a good bound.
To estimate~$\|\mu\|_\infty$ a calculation along the lines of the
bound \eqref{eq:diameter} shows that for $\mu\in C^s(D)$ with $J_n$
as in Proposition~\ref{prop:CSnonadaptive}
\begin{align}\label{eq:EstSupMu}
 \|\hat\mu_{J_n}-\mu\|_{L^\infty([a,b])}=\mathcal O_P\Big(\Big(\frac{\log
n}{n}\Big)^{s/(2s+1)}u_n\Big),
\end{align}
where $u_n=w_{J_n}/\sqrt{J_n}$. Provided the supremum of $\mu$ is attained in $[a,b]$ or $\mu$ admits some positive global H\"older regularity, we conclude that $\|\mu\|_\infty$ can be
estimated by $\|\hat\mu_{J_n}\|_\infty$ with the above rate and is in particular
a consistent estimator, which is all that is needed.
For the estimation of $\rho$ we observe that it is the second largest
eigenvalue of the transition operator $P_\Delta$. 
\citet{gobetEtAl2004} estimate this eigenvalue 
in a reflected diffusion model by constructing first an
empirical transition matrix for the transition operator restricted to a finite
dimensional space and then taking the second largest eigenvalue of the
empirical transition matrix as an estimator for $\rho$, there denoted by
$\kappa_1$. They give a rate for their estimator, in particular the estimator is
consistent.

Let us finally note that the estimation of $\zeta_\alpha$ can be circumvented by a Bayesian approach as studied by \citet{castilloNickl2014} as well as \citet{SzabovanderVaartvanZanten2014} in simpler
statistical problems.
The papers analyse Bayesian credible sets in the density estimation model and in the white noise regression model as well as in the Gaussian sequence model
and show that they are frequentist confidence sets. 
Estimating the drift of a diffusion from low-frequency observations is a
more complicated statistical model. Consistency of the
Bayesian approach in this setting has been established by
\citet{vanderMeulenvanZanten2013} and has been extended to the multi-dimensional case by \citet{GugushviliSpreij2014}.
Recently \citet{NicklSoehl2015} have shown Bayesian posterior contraction rates for scalar diffusions with unknown drift and unknown diffusion coefficient observed at low frequency.

\section{Application to diffusion processes}\label{sec:AppDiff}

\subsection{Estimation of the invariant density and its consequences}

We now apply the results from the previous section to diffusion
processes. At the same time we extend the results from inference on the
invariant probability measure to confidence bands for the drift function.
Let us consider the diffusion 
\begin{equation}
\d X_{t}=b(X_{t})\d t+\sigma \d W_{t},\qquad t\ge0,\quad X_0=x,\label{eq:sde}
\end{equation}
with a Brownian motion $W_{t}$, an unknown drift function $b\colon\R\to\R$, a volatility parameter $\sigma>0$ and starting point $x\in\R$. We observe $X$ at equidistant time
points $0,\Delta,2\Delta,\dots,(n-1)\Delta$ for some fixed observation
distance $\Delta>0$ and sample size $n\to\infty$. Our aim is inference
on the drift~$b$.

Underlying the sequence of observations $(X_{\Delta k})_{k\ge0}$ is 
a Markov structure described by the transition operator
\[
  P_{\Delta}f(x):=\E[f(X_{\Delta})|X_{0}=x].
\]
The semi-group $(P_{t}:t\ge0)$ has the infinitesimal generator $L$
on the space of twice continuously differentiable functions given
by 
\begin{equation}
Lf(x)=L_{b}f(x):=b(x)f'(x)+\tfrac{\sigma^2}{2}f''(x).\label{eq:generator}
\end{equation}
If there is an invariant density $\mu=\mu_{b}$, the operator $L$ is symmetric with respect to the scalar product of $L^{2}(\mu)=\{f:\int|f|^{2}\d\mu<\infty\}$. We impose the following assumptions on the diffusion:
\begin{assumption}\label{assDiff}
  In model~(\ref{eq:sde}) let $b$ be continuously differentiable and satisfy 
$b\in C^s(D)$ for $s\ge 1$ and a sufficiently large set $D\subset\R$ containing 
the interval $[a,b]$ for $a<b$. Let $\sigma$ be in a fixed bounded interval away from the origin. Suppose that $b'$ is bounded and that there are 
$M,r>0$ such that
  \begin{align*}
    \sgn(x)b(x)&\le-r,\quad\text{for all } |x|\ge M.
  \end{align*}
\end{assumption}
More precisely, we will need $D=[a-2^{1-j_0}(2N-1),b+2^{1-j_0}(2N-1)]$. Due to 
the global Lipschitz continuity and the assumptions on the drift, 
equation~(\ref{eq:sde}) has a unique strong solution. Moreover, $X_{t}$ is a 
Markov process with invariant probability density 
given by
\begin{equation}
\mu(x)=C_{0}\sigma^{-2}\exp\Big(2\sigma^{-2}\int_{0}^{x}b(y)\d y\Big),\quad x\in\R,\label{eq:invariantMeasure}
\end{equation}
with normalization constant $C_{0}>0$, cf. \citet[Chaps.~1,4]{bass1998}.
The corresponding Markov chain~$Z$ with $Z_k=X_{k\Delta}$ satisfies Assumption~\ref{ass:Markovchain} from the previous section.

\begin{prpstn}\label{prop:check}
If the diffusion process \eqref{eq:sde} satisfies Assumption~\ref{assDiff}, then the Markov chain $(X_{k\Delta})_{k\ge 0}$ satisfies Assumption~\ref{ass:Markovchain} where $\mu\in C^{s+1}(D)$.
\end{prpstn}
\begin{proof}
By a time-change argument we can set $\sigma=1$ without loss of generality.
\citet[Thm. 13.2]{gihmanSkorohod1979} have given an explicit formula for the transition density $p_\Delta(x,y)$ with respect to the Lebesgue measure, i.e., $P_\Delta(x,B)=\int _Bp_\Delta(x,y)\d y$ for all $B\in\mathcal B(\R)$. In particular, $p_\Delta(x,y)$ is strictly positive and thus $Z$ is $\psi$-irreducible, where $\psi$ is given by the Lebesgue measure on $\R$.

Moreover, $(x,y)\mapsto p_\Delta(x,y)$ is continuous such that for any compact interval $C\subset \R$ we have $\delta:=\delta(C):=\inf_{x,y\in C}p_\Delta(x,y)>0$
and the small set condition~\eqref{con:smallset} is satisfied:
\[
  P_\Delta(x,B)=\int_Bp_\Delta(x,y)\d y\ge \delta\int_{B\cap C}\d y=\delta|C|\nu(B),
\]
where $|C|$ denotes the Lebesgue measure of $C$ and $\nu$ is the uniform 
distribution on $C$. It also follows that the Markov chain is strongly 
aperiodic.

To show the drift condition~\eqref{con:drift}, we first construct a Lyapunov 
function for the infinitesimal generator (which is the continuous time analogue 
of the drift operator $P-\Id$), that is we find a function $V\ge1$ such that
\begin{equation}\label{eq:lyapunov}
  LV(x)\le -\lambda V(x)+c\1_C(x),\quad x\in\R.
\end{equation}
Let $V$ be a smooth function with $V(x)=e^{a|x|}$ for $|x|>R$ for some $R>0$. Due to the assumptions on $b$, we then obtain for these $x$ and $R$ large enough
\[
  LV(x)=\frac12V''(x)+b(x)V'(x)=\left(\frac{a^2}2+a\sgn (x)b(x)\right)V(x)\le
-\lambda V(x)
\]
for sufficiently small $a,\lambda$ and thus the previous inequality is 
satisfied with $C=[-R,R]$. To carry this result over to the drift 
condition~\eqref{con:drift}, we adopt the approach by \citet[Prop. 
6.4]{galtchoukPergamenchtchikov2012}: It\^o's formula yields for all $0\le 
t\le\Delta$
\[
  V(X_t)=V(x)+\int_0^t L(V)(X_s)\d s+\int _0^t V'(X_s)\d W_s.
\]
We note that Fubini's theorem yields $\E_\mu[\int_0^\Delta V'(X_s)^2\d
s]=\int_0^\Delta\E_\mu[V'(X_0)^2]\d s<\infty$ for constants~$a$ small enough by
\eqref{eq:invariantMeasure} and by the assumptions on $b$. Consequently we have
$\E_x[\int_0^\Delta V'(X_s)^2\d s]<\infty$ for almost all $x\in\R$. By the
explicit formula of $p_\Delta(x,y)$ we conclude that $\E_x [\int_0^\Delta V'(X_s)^2\d
s]<\infty$ for all $x\in\R$. Hence, the stochastic integral is a martingale (under $P_x$) and 
$Z(t):=P_t V(x)$ satisfies
\[
  Z'(t)=\E_x[L(V)(X_t)]=-\lambda Z(t)+\psi(t),\quad \psi(t):=\E_x[L(V)(X_t)+\lambda V(X_t)],
\]
where we have $\psi(t)\le cP_x(X_t\in C)\le c$ by \eqref{eq:lyapunov}. Solving this differential equation, we obtain for all $t\in[0,\Delta]$
\begin{align*}
  Z(t)=Z(0)e^{-\lambda t}+\int_0^t e^{-\lambda(t-s)}\psi(s)\d s
      \le V(x)e^{-\lambda t}+c\frac{1-e^{-\lambda\Delta}}{\lambda}.
\end{align*}
Therefore, the drift condition follows:
\[
  P_\Delta V(x)-V(x)\le (e^{-\lambda\Delta}-1)V(x)+\frac{c}{\lambda}
    \le-\tilde\lambda V(x)+\frac{c}{\lambda}\1_{\{|x|\le R\}}(x),
\]
where $R>0$ and $\tilde\lambda>0$ are chosen such that 
$(1-e^{-\lambda\Delta}-\tilde\lambda)V(x)>c/\lambda$ for $|x|>R$. In 
combination with the $\psi$-irreducibility the drift condition shows that the 
Markov chain is positive Harris recurrent.

Since our diffusion is symmetric, in the sense that the transition operator is
symmetric with respect to $L^2(\mu)$, we argue as \citet[Sect.
4.3]{bakryEtAl2008}, using that the Poincar\'e inequality is implied by
a Lyapunov--Poincar\'e inequality and we thus have the contraction 
property~\eqref{con:contraction} \citep[Thm. 1.3]{bakryEtAl2008}. Finally, the 
smoothness
of $b$ in combination with the formula for the invariant probability
density~\eqref{eq:invariantMeasure} imply that $\mu$ is in $C^{s+1}(D)$.
\end{proof}

Theorem~\ref{thrm:ucltMu} and Proposition~\ref{prop:CSnonadaptive} yield immediately
\begin{crllr}\label{cor:CSmu}
  Grant Assumption~\ref{assDiff} and let $w=(w_j)$ be increasing and satisfy $\sqrt j/w_j\to0$. Then the wavelet projection estimator $\hat\mu_{J_n}$ from \eqref{eq:hatMu} with $2^{J_n}=(n/\log n)^{1/(2s+3)}$ satisfies
  \[
    \sqrt n(\hat\mu_{J_n}-\mu)\overset{d}{\longrightarrow}\mathbb G_\mu\quad\text{in }\mathcal M_0(w).
  \]
  Moreover, the confidence band $\bar{\mathcal C}_{n}=\bar{\mathcal C}_{n}(\zeta_\alpha,s+1,u_n)$ from \eqref{eq:CS} with critical value $\zeta_\alpha$ such that $P(\|\mathbb G_\mu\|_{\mathcal M}\ge\zeta_\alpha)\le\alpha$ and $u_n=w_{J_n}/\sqrt{J_n}$ satisfies
  \[
    \liminf_{n\to\infty}P(\mu\in\bar{\mathcal C}_{n})\ge1-\alpha\quad\text{and}\quad
    |\bar{\mathcal C}_{n}|_{\infty}=\mathcal {O}_{P}\Big(\Big(\frac{n}{\log n}\Big)^{-(s+1)/(2s+3)}u_{n}\Big).
  \]
\end{crllr}

\subsection{Drift estimation via plug-in}
Supposing from now on that $\sigma=1$ and rewriting the formula of the invariant measure (\ref{eq:invariantMeasure}), we see that
\begin{equation}
b(x)=\frac{1}{2}\big(\log\mu(x)\big)'.\label{eq:drift}
\end{equation}
Obviously, $b$ depends on $\mu$ in a nonlinear way and the estimation problem is ill-posed because $b$ is a function of the derivative
$\mu'$. In general, the same calculation leads to a formula for the function $b(x)/\sigma^2$. Note that all shape properties of the drift function, like monotonicity, extrema, etc. are already determined by $b/\sigma^2$. As demonstrated by \citet{gobetEtAl2004}, the information on $\sigma$ is encoded in the transition operator of the underlying Markov chain. However, the estimation procedure in this latter article is quite involved and the construction of adaptive confidence bands in the general setting is beyond the scope of the present article. In the following we always set $\sigma=1$. Note that if we have an estimator for $\sigma$ at hand, for instance from a short high-frequency time series of the diffusion, the results easily carry over to an unknown volatility $\sigma>0$. 

Denoting the set of continuous functions on the real line by $C(\R)$, we introduce the map
$$\xi\colon \big\{f\in C^1(\R):f>0,\|f\|_{L^1}=1\big\}\to C(\R),\quad f\mapsto \frac{f'}{2f},$$
which is one-to-one with inverse function $\xi^{-1}(g)=\exp(2\int_0^\cdot g(y)\d y-c_g)$ with normalization constant $c_g\in\R$ and for any function $g$ in the range of $\xi$. We can thus estimate the drift function of the diffusion by the plug-in estimator
$\xi(\hat\mu_{J_n})$. 

Using the confidence set $\bar{\mathcal C}_{n}(\zeta_\alpha,s+1,u_n)$ for the invariant density $\mu$ from \eqref{eq:CS}, a confidence band for the drift can be constructed via
\begin{equation}\label{eq:ConfSetD}
  {\mathcal D}_{n}:={\mathcal D}_{n}(\zeta_\alpha,s,u_n):=\big\{\xi(f):f\in\bar{\mathcal C}_{n}(\zeta_\alpha,s+1,u_n)\big\}.
\end{equation}
Since $\xi$ is one-to-one, an immediate consequence of Corollary~\ref{cor:CSmu} is that we have for the coverage probability 
$\liminf_{n\to\infty}P(b\in {\mathcal D}_{n})= \liminf_{n\to\infty}P(\mu\in \bar{\mathcal C}_{n})\ge 1-\alpha.$ To 
bound the diameter of ${\mathcal D}_{n}$, we first note that $\xi$ is locally 
Lipschitz continuous: For $f,g\in C^1(\R)$ both bounded away from zero on 
$[a,b]$ we have in $L^\infty([a,b])$
\begin{align}
    \|\xi(f)-\xi(g)\|_{\infty}
    &=\frac{1}{2}\Big\|\frac{f'}{f}-\frac{g'}{g}\Big\|_{\infty}
      \le 
\frac{1}{2}\Big\|\frac{f'-g'}{f}\Big\|_{\infty}+\frac{1}{2}\Big\|\frac{
g'}{fg}(g-f)\Big\|_{\infty}\notag\\
    &\le 
\frac{1}{2}\|f^{-1}\|_{\infty}\|f'-g'\|_{\infty}
+\|\xi(g)\|_{\infty}\|f^{-1 
}\|_{\infty}\|f-g\|_{\infty}\notag\\
    &\le 
\|f^{-1}\|_{\infty}\big(\tfrac{1}{2}+\|\xi(g)\|_{\infty}
\big)\|f-g\|_{C^1([a,b])} .\label{eq:xiLipschitz}
\end{align}
For $f,g\in \bar{\mathcal C}_{n}$ we conclude in $L^\infty([a,b])$
\begin{align*}
    \|\xi(f)-\xi(g)\|_{\infty}
    &\le 
\|\xi(f)-\xi(\mu)\|_{\infty}+\|\xi(g)-\xi(\mu)\|_{\infty}\\
    &\le
\big(\tfrac{1}{2}+\|b\|_{\infty}\big)\big(\|f^{-1}
\|_{\infty}\|f-\mu\|_{C^1([a,b])}+\|g^{-1}\|_{\infty}\|g-\mu\|_{C^1([a,b])}
\big).
\end{align*}
Analogously to \eqref{eq:diameter} the choice $2^{J_n}=(n/\log n)^{1/(2s+3)}$ yields 
\[\|f-\mu\|_{C^1([a,b])}=\mathcal O_P\Big(\Big(\frac{n}{\log
n}\Big)^{-s/(2s+3)}u_n\Big)\quad\text{for all }f\in \bar
C_n(\zeta_\alpha,s+1,u_n).
\]
We conclude that $f^{-1}$ is uniformly bounded in $L^\infty([a,b])$ for all
$f\in\bar{\mathcal C}_{n}$.
Hence, we have proved
\begin{prpstn}\label{prop:onetoone}
  Grant Assumption~\ref{assDiff} with $\sigma=1,s>0$ and let $w=(w_j)$ satisfy $\sqrt 
j/w_j\to0$. Then the confidence set ${\mathcal D}_{n}={\mathcal 
D}_{n}(\zeta_\alpha,s,u_n)$ from \eqref{eq:ConfSetD} with critical value 
$\zeta_\alpha$ satisfying $P(\|\mathbb G_\mu\|_{\mathcal 
M}\ge\zeta_\alpha)\le\alpha$, $u_n=w_{J_n}/\sqrt{J_n}$ and $J_n$ chosen such that 
$2^{J_n}=(n/\log n)^{1/(2s+3)}$ fulfils
  \[
    \liminf_{n\to\infty}P(b\in\mathcal{D}_{n})\ge1-\alpha\quad\text{and}\quad
    \big|{\mathcal D}_{n}\big|_\infty=\mathcal O_P\Big(\Big(\frac{n}{\log n}\Big)^{-s/(2s+3)}u_n\Big).
  \]  
\end{prpstn}
Let us comment on the rate appearing in the previous proposition. Since the 
identification~\eqref{eq:drift} incorporates the derivative of the invariant 
measure, drift estimation is an inverse problem, which is ill-posed of 
degree one. Therefore, the minimax rate for the pointwise or $L^2$-loss is 
$n^{-s/(2s+3)}$. Considering the uniform loss, we obtain the rate $(n/\log 
n)^{-s/(2s+3)}$. Finally, $u_n\to\infty$ is the payment for undersmoothing (by 
using a weighting sequence slightly larger than $\sqrt j$). Note that we obtain 
a faster rate than \citet{gobetEtAl2004} who have proved that the minimax rate 
for drift estimation for the mean integrated squared error is $n^{-s/(2s+5)}$ 
if there is additionally an unknown volatility function in front of the 
Brownian motion in \eqref{eq:sde}.

In fact the map $\xi$ is not only Lipschitz continuous, but even Hadamard 
differentiable (on appropriate function spaces) with derivative at $\mu$
\begin{equation}\label{eq:hadamardDerivative}
  \xi'_\mu(h)=\frac{1}{2}\left(\frac{h}{\mu}\right)',\quad h\in\mathcal M_0(w).
\end{equation}
Using the delta method \cite[Thm. 20.8]{vanderVaart1998}, we obtain a functional central limit theorem for the plug-in estimator $\xi(\hat\mu_{J_n})$. 

\begin{thrm}\label{thrm:deltaMethod}
  Grant Assumption~\ref{assDiff} with $\sigma=1$ and let $w=(w_{j})$ be increasing and
  satisfy $\sqrt{j}/w_{j}\to0$, $w_j\le2^{j\delta}$ for some
$\delta\in(0,1/2)$. Let $J_{n}\in\N$ fulfil, for some $\tau\in(0,1)$,
  \[
2^{(9/4+\delta/2)J_n}w_{J_n}n^{-1/2}=o(1),\qquad
  \sqrt{n}2^{-J_{n}(2s+3)/2}w_{J_{n}}^{-1}=o(1),\qquad(\log
n)^{2/\tau}n^{-1}2^{J_{n}}J_{n}=\mathcal{O}(1).
  \]
  For $\tilde w_j\ge 2^{j(1+\delta)}$ we have as $n\to\infty$
  \[ 
 \sqrt{n}(\xi(\hat\mu_{J_n})-b)\overset{d}{\longrightarrow}\xi'_\mu(\mathbb{G}
_ { \mu })
\quad\text{in }\mathcal{M}_{0}(\tilde w).
  \]
\end{thrm}
The proof of this theorem is postponed to Section~\ref{sec:proofDeltaMethod}. 
Similarly as in \eqref{eq:CS} confidence bands for the drift function can 
alternatively be constructed by
\begin{align*}
 \bar{\mathcal D}_{n}(\bar\zeta_\alpha,s,u_n):=\left\{f:\|f-\xi(\hat\mu_{J_n})\|_{\mathcal
M(\tilde w)}<\frac{\bar \zeta_\alpha}{\sqrt{n}},\|f\|_{C^s}\le
u_n\right\},
\end{align*}
for $\alpha\in(0,1/2)$, a quantile $\bar\zeta_\alpha$ such that $P(\|\xi'_\mu(\mathbb
G_\mu)\|_{\mathcal M(\tilde w)}<\bar\zeta_\alpha)\ge1-\alpha$ and a sequence
$u_n\to\infty$. With $2^{J_n}=(n/\log n)^{1/(2s+3)}$ and
$u_n=\tilde w_{J_n}2^{-J_n}/\sqrt{J_n}$ this leads to asymptotic coverage of at least
$1-\alpha$ and a diameter decaying at rate $(n/\log n)^{-s/(2s+3)}u_n$. Note
that
in contrast to ${\mathcal D}_{n}$ the diameter of $\bar{\mathcal D}_{n}$
is slightly suboptimal due to the $\delta>0$ that appeared in
Theorem~\ref{thrm:deltaMethod} and which presumably could be removed by a more technical proof. Based on a direct estimator of the drift we will construct a similar confidence band with the optimal diameter (up to undersmoothing) in the next section.

Comparing both constructions of confidence sets, we see that ${\mathcal D}_{n}$ can be understood
as the variance stabilised version of $\bar{\mathcal D}_{n}$: The critical value of ${\mathcal D}_{n}$
depends on the unknown $\mu$ only through the covariance structure of the limit
processes $\mathbb G_\mu$ which seems to be unavoidable due to the underlying
Markov chain structure. In contrast $\bar \zeta_\alpha$ depends additionally
on $\mu$ through the derivative $\xi'_\mu$. As a consequence the confidence band~$\bar{\mathcal D}_{n}$ has the same diameter everywhere while the diameter of~${\mathcal D}_{n}$ changes.

\subsection{A direct approach to estimate the drift}

Instead of relying only on the estimator $\hat\mu_{J_n}$ of the invariant density and the plug-in approach, we can use a direct approach to estimate the drift and to obtain its confidence bands. Although there is a one-to-one correspondence between the drift function $b$ and the invariant measure $\mu$, the drift is both the canonical parameter of our model and the main parameter of interest in the context of diffusions. Since we aim for adapting to the regularity of $b$, the direct estimation approach is natural and, additionally, the resulting confidence bands will have a constant diameter.

Motivated by formula (\ref{eq:drift}), we define
our drift estimator for integers $J,U>0$ as 
\begin{equation}
\hat{b}_{J,U}=\frac{1}{2}\pi_{J}\big(\log\hat{\mu}_{J+U}\big)'=\frac{1}{2}\sum_{j\le J}\sum_{k\in K_j}\langle\big(\log\hat{\mu}_{J+U}\big)',\psi_{j,k}\rangle\psi_{j,k},\label{eq:bHat}
\end{equation}
using the wavelet projection estimator $\hat{\mu}_{J+U}$ from (\ref{eq:hatMu}). 
In contrast to the plug-in estimator in the previous section the underlying 
bias-variance trade-off is now driven by the estimation problem of~$b$ and 
the outer projection $\pi_J$ onto level $J$. However, in order to linearise the 
estimation error, we need a stable prior estimator of $\mu$ such that we
cannot simply use the empirical measure $\mu_n$ but instead use its projection 
onto some resolution level $J+U$ which is strictly
larger than $J$. As a rule of thumb, $U=U_{n}$ can be
chosen such that $2^{U_n}=\log n$ implying that an additional bias term from estimating $\mu$
is negligible. Linearising the estimation error, we obtain
\begin{equation}
\langle\hat{b}_{J,U}-b,\psi_{j,k}\rangle=-\big\langle\hat{\mu}_{J+U}-\mu,\frac{\psi_{j,k}'}{2\mu}\big\rangle+\langle R_{J+U},\psi_{j,k}\rangle,\quad j\le J,k\in K_j,\label{eq:mainStochErr}
\end{equation}
where the remainder is of order $o_{P}(n^{-1/2})$ for appropriate
choices of $J=J_{n}$, cf. Lemma~\ref{lem:linearisation} below.
In view of the linear error term and our findings in Section~\ref{sec:ucltMu},
the limit process $\mathbb{G}_{b}$ in the multi-scale space $\mathcal{M}_{0}$
will be given by 
\begin{equation}
\mathbb{G}_{b}({j,k})\sim\mathcal{N}(0,\Sigma_{j,k})\quad\text{where
}\Sigma_{j,k}:=\Sigma_{f,f}\text{ is given by \eqref{eq:Variance} with
}f=\psi_{j,k}'/(2\mu)\label{eq:limitB}
\end{equation}
with covariances $\E[\mathbb{G}_{b}({j,k})\mathbb{G}_{b}({l,m})]=\Sigma_{f_1,f_2}$ from \eqref{eq:Variance} with $f_1=\psi_{j,k}'/(2\mu)$ and $f_2=\psi_{l,m}'/(2\mu)$.
The ill-posedness of the problem is reflected by $\psi_{j,k}'$
being a factor $2^{j}$ larger than $\psi_{j,k}$. We thus need larger weights 
for high resolution levels to ensure that $\mathbb{G}_{b}$ takes values in 
$\mathcal{M}_{0}(w)$.
\begin{dfntn}\label{def:AdmissibleOne}
  A weighting sequence $w=(w_{j})$ is called
  \emph{admissible}, if it is monotonously increasing, satisfies $\sqrt{j}2^{j}/w_{j}\to0$
  as $j\to\infty$ and if there is some $\delta\in(1,2]$ such that
  $j\mapsto2^{j\delta}/w_{j}$ is monotonously increasing for large $j$.
\end{dfntn}
The last condition in the definition is a mild technical assumption
that we will need in the multi-scale central limit theorem below.
For instance, any weighting sequence $w_{j}=u_{j}\sqrt{j}2^{j}$ with
$u_{j}=j^{p}$ for some polynomial rate $p>0$ is admissible of degree
one for any $\delta\in(1,2]$. Note that admissibility of $w$ implies in
particular that $w_{j}\lesssim2^{j\delta}$ which allows to compare
the $\|\cdot\|_{\infty}$-norm with the $\|\cdot\|_{\mathcal{M}}$-norm.
We find an analogous result to Lemma~\ref{lem:PropLimit}, cf. \citet[Prop. 3]{castilloNickl2014}.
\begin{lmm}\label{lem:PropLimitB}
  $\mathbb{G}_{b}$ from (\ref{eq:limitB}) satisfies
  $\E[\|\mathbb{G}_{b}\|_{\mathcal{M}(w)}]<\infty$ for the weights
  $w$ given by $w_{j}=\sqrt{j}2^{j}$. Moreover, $\mathcal{L}(\mathbb{G}_{b})$
  is a tight Gaussian Borel probability measure in $\mathcal{M}_{0}(w)$
  for any admissible sequence $w$. 
\end{lmm}
For the following result suppose that the wavelet basis $(\phi_{j_0,l},\psi_{j,k}:j\ge j_0,l\in L,k\in K_j)$ of $L^2(\R)$ is sufficiently regular (i.e., satisfies (\ref{eq:fPsi}) with $\gamma\ge3/2+\delta$), for instance, Daubechies' wavelets of order $N\ge20$.
\begin{thrm}\label{thrm:ucltDrift}
  Grant Assumption~\ref{assDiff} with $\sigma=1$ and let $w=(w_{j})$ be admissible. Let $J_{n}\to\infty,U_{n}\to\infty$ fulfil
  \[
  \sqrt{n}2^{-J_{n}(s+1/2)}w_{J_{n}}^{-1}=o(1),\quad\sqrt{n}2^{-(J_{n}+U_{n})(2s+1)}=o(1),\quad n^{-1/2}2^{2(J_{n}+U_{n})}(J_{n}+U_{n})=o(1).
  \]
  Then $\hat{b}_{J,U}$ from (\ref{eq:bHat}) satisfies, as $n\to\infty$,
  \[
  \sqrt{n}(\hat{b}_{J_{n},U_{n}}-b)\overset{d}{\longrightarrow}\mathbb{G}_{b}\quad\text{in }\mathcal{M}_{0}(w)
  \]
  for the tight Gaussian random variable $\mathbb{G}_{b}$ in $\mathcal{M}_{0}(w)$
  given by (\ref{eq:limitB}).
\end{thrm}
The proof of this theorem is postponed to Section~\ref{sec:ProofUCLT}.
The first condition on $J_{n}$ is the bias condition for $b$ in
$\mathcal{M}_{0}.$ The latter two conditions on $J_{n}+U_{n}$ are
determined by a bias and a variance condition for $\mu$ which we will
need to bound the remainder $R_{J_{n}+U_{n}}$ from (\ref{eq:mainStochErr})
in $L^{\infty}$. If $\delta<1/2+s$ in Definition~\ref{def:AdmissibleOne},
then the second condition is strictly weaker than the first one.

Similarly to the confidence band for $\mu$ in
Proposition~\ref{prop:CSnonadaptive} we can now construct a confidence band for
the drift function $b$.
For some $\alpha\in(0,1)$ we consider
\begin{equation}\label{eq:CSb}
  \mathcal E_n:=\mathcal E_{n}(\zeta_{\alpha},s,u_n):=\Big\{ f:\|f-\hat{b}_{J_{n},U_{n}}\|_{\mathcal M}<\frac{\zeta_{\alpha}}{\sqrt{n}}, \|f\|_{C^{s}}\le u_{n}\Big\},
\end{equation}
where $\zeta_{\alpha}$ is chosen such that
$P(\|\mathbb{G}_{b}\|_{\mathcal{M}}<\zeta_{\alpha})\ge1-\alpha$ and $(u_{n})_n$ is a diverging sequence.
\begin{prpstn}\label{prop:CSnonadaptiveb}
  Grant Assumption~\ref{assDiff} with
  $\sigma=1,s\ge1$ and let $w=(w_{j})$ be admissible. For $\alpha\in(0,1)$
  let $\zeta_{\alpha}>0$ satisfy $P(\|\mathbb{G}_{b}\|_{\mathcal{M}}\ge\zeta_{\alpha})\le\alpha$
  and choose $J_{n}:=J_{n}(s)$ and $U_n\to\infty$ such that 
  \[
  2^{J_{n}}=\Big(\frac{n}{\log n}\Big)^{1/(2s+3)}\quad\text{and}\quad2^{U_{n}}=\mathcal{O}(\log n).
  \]
  Then the confidence set ${\mathcal E}_{n}={\mathcal E}_{n}(\zeta_{\alpha},s,u_{n})$
  from (\ref{eq:CSb}) with $u_{n}:=w_{J_{n}}2^{-J_{n}}/\sqrt{J_{n}}$
  satisfies
  \[
  \liminf_{n\to\infty}P(b\in{\mathcal E}_{n})\ge1-\alpha\quad\text{and}\quad|{\mathcal E}_{n}|_{\infty}=\mathcal{O}_{P}\Big(\Big(\frac{n}{\log n}\Big)^{-s/(2s+3)}u_{n}\Big).
  \]
\end{prpstn}
\begin{proof}
The proof is essentially the same as for the confidence band of the invariant
probability density. We show that the asymptotic coverage probability is at least
$1-\alpha$ and obtain for $f,g\in\mathcal E_n$ as in~\eqref{eq:diameter} the bound
\begin{align}
\|f-g\|_{\infty}&
= \mathcal{O}_{P}\Big(n^{-1/2}2^{J_{n}/2}w_{J_{n}}\Big)+\mathcal{O}_{P}\Big(2^{-J_
{n}s}u_{n}\Big)\nonumber 
\end{align}
Using $u_{n}=w_{J_{n}}2^{-J_{n}}/\sqrt{J_{n}}$ we thus have
\begin{align*}
\|f-g\|_{\infty}&
= 
\mathcal{O}_{P}\Big(n^{-1/2}2^{3J_{n}/2}J_{n}^{1/2}u_{n}\Big)+\mathcal{O}_{P}
\Big(2^{-J_{n}s}u_{n}\Big).
\end{align*}
The choice of $J_{n}$ yields $n^{-1/2}\sqrt{2^{3J_{n}}J_{n}}\lesssim(n/\log
n)^{-s/(2s+3)}=2^{-J_{n}s}$.
\end{proof}

\section{Adaptive confidence bands for drift estimation}\label{sec:AdaptiveConf}

Inspired by \citet{gineNickl2010}, we will now construct an adaptive
version of the confidence set ${\mathcal E_n}$ from \eqref{eq:CSb}.
To this end we estimate the regularity $s$ of the drift with a Lepski-type 
method.
For some maximal regularity $r>1$, let the integers $0<J_{min}<J_{max}$
be given by
\[
2^{J_{min}}\sim\Big(\frac{n}{\log n}\Big)^{1/(2r+3)},\quad2^{J_{max}}\sim\frac{n^{1/4}}{(\log n)^{2}}.
\]
Note that $J_{min},J_{max}$ depend on the sample size $n$, which
is suppressed in the notation. If we knew in advance that $b$ has
regularity $r$, then we would choose the resolution level $J_{min}$.
The upper bound $J_{max}$ is chosen such that $J_{max}+U_{n}$ satisfies
the third condition in Theorem~\ref{thrm:ucltDrift}. The set in which
we will adaptively choose the optimal resolution level for regularities
$s\in[1,r]$ is defined by
\[
\mathcal{J}_{n}:=[J_{min},J_{max}]\cap\mathbb{N}.
\]
Similar to \citet[Lem. 2]{gineNickl2010}, we show under the following
assumption on $b$ that the optimal truncation level
can be consistently estimated up to a fixed integer.
\begin{assumption}
  \label{ass:SelfSim}Let $b\in C^{s}(D),s\ge1,$ satisfy for constants
  $0<d_{1}<d_{2}<\infty$ and an integer $J_{0}>0$ that
  \begin{equation}\label{eq:SelfSim}
  d_{1}2^{-Js}\le\|\pi_{J}(b)-b\|_{L^\infty([a,b])}\le d_{2}2^{-Js},\quad\forall J\ge J_{0}.
  \end{equation}
\end{assumption}
The second inequality in \eqref{eq:SelfSim} is the well known Jackson inequality 
which is satisfied for all usual choices of wavelet basis. The first inequality 
is the main condition here, called \emph{self-similarity assumption}. It 
excludes the cases where the bias would be smaller than the usual order $2^{-Js}$. 
Although the estimator $\hat b_{J,U}$ would profit from a smaller bias, we 
cannot hope for a consistent estimation of the optimal projection level and the 
resulting regularity index $s$ if \eqref{eq:SelfSim} (or a slight generalisation 
by \citet{bull2012}) is violated. Indeed, \citet{hoffmannNickl2011} have shown 
that this kind of condition is necessary to construct adaptive and honest 
confidence bands. On the other hand, it has been proved by \citet{gineNickl2010} 
that the set of functions that do not satisfy the self-similarity assumption is 
nowhere dense in the H\"older norm topologies. 
In that sense, the self-similarity assumption is satisfied by ``typical'' functions. We will give an illustrative example next. Probabilistic examples for self-similar functions are those Gaussian processes which can be represented as stochastic series expansions like the Karhunen-Lo\`eve expansion for Brownian motion or typical examples of Bayesian priors. Naturally, more regular functions $b\in C^{r}(D)$ for some $r>s$ cannot satisfy Assumption~\ref{ass:SelfSim}. For a further discussion and examples we refer to \citet[Section~3.5]{gineNickl2010} as well as \citet{bull2012}.  
\begin{xmpl}
 Let $b$ be a smooth function on $\R$ except for some point $x_0\in(a,b)$ where the $s$th order derivative $b^{(s)}$ has a jump for some integer $s\ge1$. Locally around $x_0$ the function $b$ can be approximated by a Taylor polynomial $b(x)=\beta_s^{\pm}(x-x_0)^s+\sum_{l=0}^{s-1} \beta_l(x-x_0)^l+\mathcal O(|x-x_0|^{s+1})$ with coefficients $\beta_0,\dots,\beta_{s-1}\in\R$ and where the coefficient of order $s$ is some $\beta_s^+\in\R$ or $\beta_s^-\in\R$ depending on $x\ge x_0$ and $x< x_0$, respectively. Due to the jump of $b^{(s)}$ we have $\beta_s^+\neq\beta_s^-$.  
 
 Choose $k_j$ as the nearest integer of $2^{j}x_0$, implying that $x_0$ is in the middle of the support of $\psi_{j,k_j}$. Using the elementary estimate $|\langle f,\psi_{j,k}\rangle|\le\|f\|_{L^\infty(\supp\psi_{j,k})}2^{-j/2}\|\psi\|_{L^1}$, we obtain
 \begin{align*}
    \|\pi_{J}(b)-b\|_{L^\infty([a,b])}
    &=\Big\|\sum_{j>J,k}\langle b,\psi_{j,k}\rangle\psi_{j,k}\Big\|_{L^\infty([a,b])}\gtrsim \sup_{j> J}2^{j/2}|\langle b,\psi_{j,k_j}\rangle|.
 \end{align*}
For sufficiently large $j$ the regularity of the wavelet basis, being thus orthogonal to polynomials, and their compact support yield
\begin{align*}
   2^{j/2}|\langle b,\psi_{j,k_j}\rangle|
   &\ge 2^{j}\Big|\int \big(b(x)-\beta_s^{-}(x-x_0)^s+\sum_{l=0}^{s-1} \beta_l(x-x_0)^l\big)\psi(2^j x-k_j)\d x\Big|+\mathcal O(2^{-j(s+1)})\\
   &\ge 2^j|\beta_s^+-\beta_s^-|\Big|\int_{x>x_0}(x-x_0)^s\psi(2^jx-k_j)\d x\Big|+\mathcal O(2^{-j(s+1)})\\
   &\ge 2^{-js}|\beta_s^+-\beta_s^-|\min_{\epsilon\in[-1/2,1/2]}\Big|\int_{y>0}y^s\psi(y+\epsilon)\d y\Big|+\mathcal O(2^{-j(s+1)}).
\end{align*}
We conclude for $J$ sufficiently large that $\|\pi_{J}(b)-b\|_{L^\infty([a,b])}\gtrsim 2^{-Js}$.
\end{xmpl}

The oracle choice $J_{n}^{*}$ which balances the bias $\|\pi_{J}(b)-b\|_{L^\infty([a,b])}$
and the main stochastic error is given by 
\[
J_{n}^{*}:=J_{n}^{*}(s)=\min\Big\{ J\in\mathcal{J}_{n}:(d_{2}+1)2^{-Js}\le\frac{K}{4}\sqrt{\frac{2^{3J}J}{n}}\Big\}
\]
for some suitable constant $K>0$ depending only on $\psi$, $\inf\{\mu(x):x\in\cup_{l\in L}\supp\phi_{j_0,l}\}$ and the maximal asymptotic variance $\tilde\Sigma=\sup_{j,k}\Sigma_{j,k}$ where the latter two quantities can be replaced by the consistent estimators which we have discussed in Section~\ref{sec:ConstConfBand}. We see easily that $2^{J_{n}^{*}}\sim(\frac{n}{\log n})^{1/(2s+3)}$.
Following Lepski's approach, we define the estimator for $J_{n}^{*}$
by
\begin{equation}
\hat{J}_{n}=\min\Big\{ 
J\in\mathcal{J}_{n}:\|\hat{b}_{J}-\hat{b}_{j}\|_{L^\infty([a,b])}\le 
K\sqrt{\frac{2^{3j}j}{n}}\;\forall j>J,j\in\mathcal{J}_{n}\Big\}.\label{eq:Jhat}
\end{equation}

\begin{lmm}\label{lem:Lepski}
  Grant Assumptions~\ref{assDiff} and \ref{ass:SelfSim}
  for $s\in[1,r]$ with some $r>1$ and $\sigma=1$. Let $w$ be admissible. Then there are a constant $K>0$ depending only on $\psi$, $\inf\{\mu(x):x\in\cup_{l\in L}\supp\phi_{j_0,l}\}$ and the maximal asymptotic variance $\tilde\Sigma=\sup_{j,k}\Sigma_{j,k}$, an integer $M>0$ depending only on $d_{1},d_{2},K$ and for any $\tau\in(0,1)$ there are constants $C,c>0$ depending on
  $\tau,K,\psi$ such that
  \[
  P\big(\hat{J}_{n}\notin[J_{n}^{*}-M,J_{n}^{*}]\big)\le C\big(n^{-cJ_{min}^{\tau}}+e^{-cJ_{min}}\big)\to0.
  \]
\end{lmm}
The proof of this lemma relies on the concentration result in Proposition~\ref{prop:Concentration}
and is postponed to Section~\ref{sec:ProofofLepski}. Applying that
$\hat{J}_{n}$ is a reasonable estimator of $J_{n}^{*}$, we obtain
a completely data-driven estimator 
\begin{equation}
\hat{b}:=\hat{b}_{\hat{J}_{n},U_{n}},\quad\text{with }\hat{J}_{n}\text{ from \eqref{eq:Jhat} and }2^{U_{n}}=\log n.\label{eq:adaptiveEst}
\end{equation}

\begin{crllr}\label{cor:adaptiveUCLT}
In the situation of Lemma~\ref{lem:Lepski} the adaptive estimator
$\hat{b}$ defined by (\ref{eq:adaptiveEst}) satisfies $\|\hat{b}-b\|_{\mathcal{M}}=\mathcal{O}_{P}(n^{-1/2})$
and $\|\hat{b}-b\|_{L^\infty([a,b])}=\mathcal{O}_{P}\big((n/\log n)^{-s/(2s+3)}u_{n}\big)$
with $u_{n}:=w_{J_{n}^{*}}2^{-J_{n}^{*}}/\sqrt{J_{n}^{*}}$.
Further for every $m\in\{0,1,\dots,M\}$ we have 
\[  \sqrt{n}(\hat{b}_{{J}_{n}^*-m,U_{n}}-b)\overset{d}{\longrightarrow}\mathbb{G}_{b}\quad\text{in }\mathcal{M}_{0}(w)
  \]
  as $n\to\infty$
  for the tight Gaussian random variable $\mathbb{G}_{b}$ in $\mathcal{M}_{0}(w)$
  given by (\ref{eq:limitB}).
\end{crllr}
\begin{proof}
Combining Lemma~\ref{lem:Lepski} and Theorem~\ref{thrm:ucltDrift}, there is for any $\delta>0$ a constant $C>0$ such that for $n$ large enough
\[
P\big(\sqrt{n}\|\hat{b}-b\|_{\mathcal{M}}>C\big)\le\sum_{J=J_{n}^{*}-M}^{J_{n}^{*}}P\big(\sqrt{n}\|\hat{b}_{J}-b\|_{\mathcal{M}}>C\big)+o(1)\le (M+1)\delta.
\]
Since $M$ is a finite constant, we have$\|\hat{b}-b\|_{\mathcal{M}}=\mathcal{O}_{P}(n^{-1/2})$.
Using that $2^{\hat{J}_{n}}\sim2^{J_{n}^{*}}$, a calculation similar
to (\ref{eq:diameter}) yields the bound for the uniform norm.
For the second claim notice that the estimators $\hat{b}_{{J}_{n}^{*}-m,U_{n}}$ satisfy the conditions of Theorem~\ref{thrm:ucltDrift}.
\end{proof}
The bound for the uniform risk is slightly suboptimal because $u_{n}$
diverges arbitrary slowly (depending on the choice of $w$) to infinity.
Using direct estimates of the $\|\cdot\|_{\infty}$-norm in the proofs
in Section~\ref{sec:ProofUCLT}, this additional factor could be
circumvented. However, it can be interpreted as an additional factor that corresponds to a slight undersmoothing which is often used to have a negligible bias in the construction of confidence bands.

Another consequence of Lemma~\ref{lem:Lepski} is that we can consistently
estimate the regularity~$s$ of~$b$. For a sequence of random variables $(v_n)$ with $v_{n}^{-1}=o_P(1)$
we define the estimator
\begin{equation}
\hat{s}_{n}:=\max\bigg(1,\,\frac{\log n-\log\log n}{2(\log2)(\hat{J}_{n}+v_{n})}-\frac{3}{2}\Big(1+\frac{v_{n}}{\hat{J}_{n}}\Big)\bigg).\label{eq:sHat}
\end{equation}
Using that $2^{J_{n}^{*}}\sim(n/\log n)^{1/(2s+3)}$, we derive
from Lemma~\ref{lem:Lepski} the following corollary. The proof can
be found in Section~\ref{sec:ProofSHat}.
\begin{crllr}
\label{cor:sHat}In the situation of Lemma~\ref{lem:Lepski} the
estimator $\hat{s}_{n}$ given by (\ref{eq:sHat}) satisfies  for any sequence of random variables 
$(v_{n})$ with $v_{n}^{-1}=o_P(1)$
\[
P(\hat{s}_{n}\le s)\to1\quad\text{and}\quad s-\hat{s}_{n}=\mathcal{O}_{P}\Big(\frac{v_{n}}{J_{n}^{*}}\Big).
\]
\end{crllr}
With the estimator $\hat{b}$ from (\ref{eq:adaptiveEst}) we can now construct
our adaptive confidence bands as follows. 
By a Bonferroni correction we take care of the possible dependence between the estimators $\hat{b}_{{J}_{n}^*-m,U_{n}}$ and the adaptive choice $\hat J_n$ of the resolution level. 
In this way sample splitting can be avoided, which was also used by \citet{bull2012}.
For any level $\alpha\in(0,1)$ let $\beta=\alpha/(M+1)$ and
define
\begin{equation}
\tilde{\mathcal E}_{n}:=\tilde{\mathcal E}_{n}(\hat{\zeta}_{\beta,n},\hat{s}_{n},t_{n}):=\Big\{ f:\|f-\hat{b}\|_{\mathcal{M}}<\frac{\hat{\zeta}_{\beta,n}}{\sqrt{n}},\|f\|_{C^{\hat{s}_{n}}}\le t_{n}\Big\},\label{eq:adaptiveCS}
\end{equation}
where $(t_{n})$ is a sequence of random variables with $t_{n}^{-1}=o_P(1)$ and $\hat{\zeta}_{\beta,n}$
is an (over-)estimator of the critical value $\zeta_{\beta}$ given by $P(\|\mathbb{G}_{b}\|_{\mathcal{M}}<\zeta_{\beta})\ge1-\beta$ similarly to the construction in Section~\ref{sec:ConstConfBand}. Now we can state our final theorem:
\begin{thrm}\label{thrm:adap}
  Grant Assumptions~\ref{assDiff} and \ref{ass:SelfSim} for $\sigma=1,s\in[1,r]$ with some $r>1$. Let $w=(w_{j})$ be admissible and define $u_{n}:=w_{J_{n}^{*}}2^{-J_{n}^{*}}/\sqrt{J_{n}^{*}}$ as well as $\hat u_n:=w_{\hat J_{n}}2^{-\hat J_{n}}\hat J_{n}^{-1/2}$. For $\alpha\in(0,1)$ set $\beta:=\alpha/(M+1)$ with $M$ as in Lemma~\ref{lem:Lepski}. Let $\zeta_{\beta}>0$ be given by $P(\|\mathbb{G}_{b}\|_{\mathcal{M}}\ge\zeta_{\beta})\le\beta$
  and let $\hat{\zeta}_{\beta,n}$ be an (over-)estimator satisfying $P(\hat{\zeta}_{\beta,n}\ge\zeta_{\beta}-\varepsilon)\to1$ for all $\varepsilon>0$.
  Then the confidence set $\tilde{\mathcal E}_{n}=\tilde{\mathcal E}_{n}(\hat{\zeta}_{\beta,n},\hat{s}_{n},t_{n})$
  given by~(\ref{eq:adaptiveCS}), where we choose
$t_{n}:=\sqrt{\hat u_n}$ 
  and $\hat{s}_{n}$
  according to~(\ref{eq:sHat}) with $v_n=o_P(\log \hat u_n)$ and $v_{n}^{-1}=o_P(1)$, satisfies
  \[
  \liminf_{n\to\infty}P(b\in\tilde{\mathcal E}_{n})\ge 1-\alpha\quad\text{and}\quad|\tilde{\mathcal E}_{n}|_{\infty}=\mathcal{O}_{P}\Big(\Big(\frac{n}{\log n}\Big)^{-s/(2s+3)}u_{n}\Big).
  \]
\end{thrm}
\begin{proof}
We will adapt the proof of Proposition~\ref{prop:CSnonadaptiveb}
to the estimated quantities $\hat{J}_{n},\hat{s}_{n}$ and $\hat{\zeta}_{\beta}$.
By Corollary~\ref{cor:sHat} the probability of the event $\{\hat{s}_{n}\le s\}$
converges to one. Due to $v_{n}^{-1}=o_P(1)$, we thus have $b\in C^{\hat{s}_{n}}(v_{n})$
with probability tending to one. Using Lemma~\ref{lem:Lepski} and Corollary~\ref{cor:adaptiveUCLT}, we infer
\begin{equation*}
\limsup_{n\to\infty}P(b\notin\tilde{\mathcal E}_{n})
= 
\limsup_{n\to\infty}\sum_{m=0}^{M}P\Big(\sqrt{n}\|\hat{b}_{J_n^*-m,U_n}-b\|_{\mathcal{M}}\ge\hat{\zeta}_{\beta,n}\Big)\le (M+1)P(\|\mathbb{G}_{b}\|_{\mathcal{M}}\ge\zeta_{\beta})\le\alpha.
\end{equation*}
We conclude that $  \liminf_{n\to\infty}P(b\in\tilde{\mathcal E}_{n})\ge 1-\alpha$.

To estimate the diameter, we proceed as in (\ref{eq:diameter}). Applying
additionally Corollary~\ref{cor:sHat}, we obtain for any $f,g\in\tilde{\mathcal E}_{n}$
\begin{align*}
\|f-g\|_{\infty}&\lesssim  \sum_{j\le J_{n}^{*}}2^{j/2}\max_{k}|\langle 
f-g,\psi_{j,k}\rangle|+\sum_{j>J_{n}^{*}}2^{j/2}\max_{k}|\langle 
f-g,\psi_{j,k}\rangle|\\
&\le  
\big(\|f-\hat{b}\|_{\mathcal{M}}+\|g-\hat{b}\|_{\mathcal{M}}\big)\sum_{j\le 
J_{n}^{*}}2^{j/2}w_{j}+\|f-g\|_{C^{\hat{s}_{n}}}\sum_{j>J_{n}^{*}}2^{-j\hat{s}_{ 
n}}\\
&=  
\mathcal{O}_{P}\big(n^{-1/2}2^{J_{n}^{*}/2}w_{J_{n}^{*}}\big)+\mathcal{O}_{P}(t_
{n}2^{-J_{n}^{*}s+\mathcal{O}_P(v_{n})})\\
&=  
\mathcal{O}_{P}\Big(n^{-1/2}2^{3J_{n}^{*}/2}(J_{n}^{*})^{1/2}u_{n}+2^{-J_{n}^{*}
s}\hat u_{n}\Big),
\end{align*}
where we have plugged in the choices of $t_{n}$, $u_n$ and $v_{n}$
and $\hat u_n\lesssim u_n$ with probability converging to one. Since $2^{J_{n}^{*}}\sim(n/\log n)^{1/(2s+3)}$,
the assertion follows.
\end{proof}

The confidence bands are constructed explicitly and this helps to verify that the confidence bands are \emph{honest}, i.e. the coverage is achieved uniformly over some set of the unknown parameter. The general philosophy being that uniformity in the assumptions leads to uniformity in the statements, the detailed derivation of honesty is tedious so that we only sketch it here. The main ingredients of the proof are the central limit theorem and the concentration inequality for Markov chains. In the original version of the concentration inequality, Theorem~9 by~\citet{adamczakBednorz2013}, the constants are given explicitly in terms of the assumptions and thus the concentration inequality is uniform in the underlying Markov chain $Z$. It is also to be expected that the central limit theorem holds uniformly in the bounded-Lipschitz metric with respect to~$Z$ although this is not explicitly contained in the statement. With these uniform ingredients a uniform version of Theorem~\ref{thrm:ucltMu} can be proved, where the convergence in distribution is again metrised in the bounded-Lipschitz metric. In combination with a uniform bound on the Lebesgue densities of $\|\mathbb G_\mu\|_{\mathcal M}$ this leads to honest confidence bands in Proposition~\ref{prop:CSnonadaptive}. Thanks to the explicit derivation of Assumption~\ref{ass:Markovchain} from Assumption~\ref{assDiff}, uniformity in the diffusion model carries over to uniformity in the Markov chain and we see that the confidence bands in Proposition~\ref{prop:onetoone} are honest. Likewise a uniform version of Theorem~\ref{thrm:ucltDrift} can be proved. Provided the random variables $\|\mathbb{G}_b\|_{\mathcal{M}}$ have uniformly bounded Lebesgue densities this uniform version entails honest and adaptive confidence bands for the drift in Theorem~\ref{thrm:adap}.

\section{Proof of Theorem~\ref{thrm:ucltDrift}\label{sec:ProofUCLT}}

In the sequel we use the notation 
\[
J^{+}=J+U
\]
for the projection level of $\hat{\mu}_{J^{+}}$ and we define 
\[S:=\bigcup_{l\in L}\supp\phi_{j_0,l}\subset[a-2^{-j_0}(2N-1),b+2^{-j_0}(2N-1)].\]
To analyse the estimation error of the wavelet coefficients $\langle\hat{b}_{J,U},\psi_{j,k}\rangle$,
we apply the following linearisation lemma:
\begin{lmm}
\label{lem:linearisation}Grant Assumption~\ref{assDiff} with $\sigma=1$. For
$j\in\{-1,j_0,\dots,J\}$ and $k\in K_j$ we have
\[
\langle\hat{b}_{J,U}-b,\psi_{j,k}\rangle=-\big\langle\hat{\mu}_{J^{+}}-\mu,\frac{\psi_{j,k}'}{2\mu}\big\rangle+\langle R_{J^{+}},\psi_{j,k}\rangle,
\]
where the remainder is given by $R_{J^{+}}=-\frac{\hat{\mu}_{J^{+}}-\mu}{2\hat{\mu}_{J^{+}}}\Big(\frac{\hat{\mu}_{J^{+}}-\mu}{\mu}\Big)'$.
If $J^{+}=J_{n}^{+}$ satisfies for some $\tau\in(0,1)$ 
\[
2^{-J_{n}^{+}(s+1)}=o(1),\qquad(\log n)^{2/\tau}n^{-1}2^{J_{n}^{+}}J_{n}^{+}=\mathcal{O}(1),
\]
then
\[
\|R_{J_{n}^{+}}\|_{L^\infty(S)}=\mathcal{O}_{P}\big(n^{-1}J_{n}^{+}2^{2J_{n}^{+}}+2^{-J_{n}^{+}(2s+1)}\big).
\]
\end{lmm}
\begin{proof}
Writing $\eta:=(\hat{\mu}_{J^+}-\mu)/\mu$, the chain rule yields 
\begin{align*}
\frac{1}{2}(\log \hat\mu_{J^+})'-b & =\frac{1}{2}\big(\log(1+\eta)\big)'=\frac{\eta'}{2(1+\eta)}=\frac{1}{2}\Big(\frac{\hat{\mu}_{J^{+}}-\mu}{\mu}\Big)'+R_{J^{+}},
\end{align*}
where the remainder is given by
\[
R_{J^{+}}:=-\frac{\eta\eta'}{2(1+\eta)}=-\frac{\hat{\mu}_{J^{+}}-\mu}{2\hat{\mu}
_{J^{+}}}\Big(\frac{\hat{\mu}_{J^{+}}-\mu}{\mu}\Big)'.
\]
Using integration by parts with vanishing boundary terms, the wavelet coefficients corresponding to the linear
term can be written as 
\[
\frac{1}{2}\langle((\hat{\mu}_{J^{+}}-\mu)\mu^{-1})',\psi_{j,k}\rangle=-\frac{1}{2}\langle\hat{\mu}_{J^{+}}-\mu,\psi_{j,k}'\mu^{-1}\rangle.
\]

Let us bound the remainder, starting with $\|\hat{\mu}_{J_{n}^{+}}'-\mu'\|_{L^\infty(S)}.$
Decomposing the uniform error into a bias and a stochastic error term,
we obtain
\begin{align*}
\|\hat{\mu}_{J_{n}^{+}}'-\mu'\|_{L^\infty(S)}&\le  \Big\|\sum_{j\le 
J_{n}^{+},k}\langle\mu_{n}-\mu,\psi_{j,k}\rangle\psi_{j,k}'\Big\|_{L^\infty(S)}
+\Big\|\sum_{j>J_{n}^{+},k}\langle\mu,\psi_{j,k}\rangle\psi_{j,k}'\Big\|_{
L^\infty(S)}\\
&=:  V_{n}+B_{n}.
\end{align*}
Using the localisation property of the wavelet function $\|\sum_{k}|\psi(\bull-k)|\|_{\infty}\lesssim1$
(which holds for $\psi'$ as well) and the regularity of $\mu\in C^{s+1}(D)$,
implying $\sup_{j,k:\supp\psi_{j,k}\cap S\neq\emptyset}2^{j(s+3/2)}|\langle\mu,\psi_{j,k}\rangle|<\infty$,
the bias can be estimated by
\[
B_{n}\lesssim\sum_{j>J_{n}^{+}}\max_{k:\supp\psi_{j,k}\cap S\neq\emptyset}|\langle\mu,\psi_{j,k}\rangle|\Big\|\sum_{k}|\psi_{j,k}'|\Big\|_{L^\infty(S)}\lesssim\sum_{j>J_{n}^{+}}2^{-js}\lesssim2^{-J_{n}^{+}s}.
\]
For the stochastic error term we obtain similarly
\[
V_{n}\lesssim\sum_{j\le J_{n}^{+}}2^{3j/2}\max_{k:\supp\psi_{j,k}\cap S\neq\emptyset}|\langle\mu_{n}-\mu,\psi_{j,k}\rangle|.
\]
The maximum of $2^{j}$ subgaussian random variables is of order $\mathcal{O}_{P}(\sqrt{j})$.
More precisely, Proposition~\ref{prop:Concentration} and the assumptions
on $J_{n}^{+}$ yield for any $j_0\le j\le J_{n}^{+}$ and $\tau\in(0,1)$,
similarly to (\ref{eq:ApplyBernstein}), 
\begin{align*}
  &P\left(\max_{k:\supp\psi_{j,k}\cap 
S\neq\emptyset}\sqrt{n}|\langle\mu_{n}-\mu,\psi_{j,k}\rangle|\ge\sqrt{j}t\right)
  \lesssim  2^{j}\exp\big(-c_{1}(\log n)j^{\tau}t^{\tau}\big)+2^{j}\exp\big(-c_{2}j(t\wedge t^{2})\big)\\
  &\qquad\qquad\le \exp\big(j(\log2-c_{1}j^{\tau-1}(\log n)t^{\tau})\big)+\exp\big(j(\log2-c_{2}(t\wedge t^{2}))\big).
\end{align*}
Using $J_{n}^{+}\lesssim\log n$, the right-hand side of the previous
display is arbitrarily small for large enough~$t$. An analogous estimate holds 
for the scaling functions~$\psi_{-1,\cdot}$. Therefore,
\[
\|\hat{\mu}_{J_{n}^{+}}'-\mu'\|_{\infty}=\mathcal{O}_{P}\Big(\sum_{j_0\le j\le 
J_{n}^{+}}2^{3j/2}\sqrt{j/n}\Big)+\mathcal{O}\big(2^{-J_{n}^{+}s}\big)=\mathcal{
O}_{P}\left(2^{3J_{n}^{+}/2}\sqrt{J_{n}^{+}/n}\right)+\mathcal{O}\big(2^{-J_{n}^
{ + } s } \big).
\]
Analogously, we have
\[
\|\hat{\mu}_{J_{n}^{+}}-\mu\|_{\infty}=\mathcal{O}_{P}\left(2^{J_{n}^{+}/2}\sqrt
{ J_ {n}^{+}/n}\right)+\mathcal{O}\big(2^{-J_{n}^{+}(s+1)}\big).
\]
Since $\mu$ is bounded away from zero on $S$, the choice of $J_{n}^{+}$
yields in particular that we have $\lim_{n\to\infty}P(\inf_{x\in 
S}\hat{\mu}_{J_{n}^{+}}(x)>\frac{1}{2}\inf_{x\in S}\mu(x))=1$.
We conclude
\begin{align*}
  \|R_{J_{n}^{+}}\|_{L^\infty(S)} 
  & =\mathcal{O}_{P}\big(\|\hat{\mu}_{J_{n}^{+}}'-\mu'\|_{L^\infty(S)}\|\hat{\mu}_{J_{n}^{+}}-\mu\|_{L^\infty(S)}+\|\hat{\mu}_{J_{n}^{+}}-\mu\|_{L^\infty(S)  }^2\big)\\
  & =\mathcal{O}_{P}\Big(\big(2^{3J_{n}^{+}/2}\sqrt{J_{n}^{+}/n}+2^{-J_{n}^{+}s}\big)\big(2^{J_{n}^{+}/2}\sqrt{J_{n}^{+}/n}+2^{-J_{n}^{+}(s+1)}\big)\Big),
\end{align*}
which shows the asserted bound for $\|R_{J_{n}^{+}}\|_{L^\infty(S)}$.
\end{proof}
The linearised stochastic error term can be decomposed into
\begin{align}
-\frac{1}{2}\langle\hat{\mu}_{J^{+}}-\mu,\psi_{j,k}'\mu^{-1}\rangle&=  -\sum_{l\le J^{+},m}\langle\mu_{n}-\mu,\psi_{l,m}\rangle\big\langle\psi_{l,m},\frac{\psi_{j,k}'}{2\mu}\big\rangle+\big\langle(\Id-\pi_{J^{+}})\mu,\frac{\psi_{j,k}'}{2\mu}\big\rangle\nonumber \\
&=  -\big\langle\mu_{n}-\mu,\frac{\psi_{j,k}'}{2\mu}\big\rangle+\big\langle\mu_{n}-\mu,(\Id-\pi_{J^{+}})\Big(\frac{\psi_{j,k}'}{2\mu}\Big)\big\rangle+\big\langle(\Id-\pi_{J^{+}})\mu,\frac{\psi_{j,k}'}{2\mu}\big\rangle.\label{eq:DecompLin2}
\end{align}
Roughly for $j\le J\le J^{+}$, Theorem~\ref{thrm:ucltMu} (or an
analogous result for ill-posed problems) applies for the first term
in the above display, the second term should converge to zero by the
localisation of the~$\psi_{j,k}$ in the Fourier domain and the third
term is a bias that can be bounded by the smoothness of $\mu$. If
$U_{n}\to\infty$ this ``$\mu$-bias'' term is of smaller order
than the ``$b$-bias'' which is determined by $\sum_{j>J,k}\langle b,\psi_{j,k}\rangle\psi_{j,k}$.

Let us make these considerations precise. We will need the following
lemma, which relies on the localisation of the wavelets in Fourier
domain. More precisely, $\psi$ can be chosen such that for some $\gamma\ge1$
we have 
\[
\phi,\psi\in C^{\gamma}(\R)\quad\text{and}\quad\int x^{k}\psi(x)\d x=0,\quad\text{for }k=0,\dots,\lceil\gamma\rceil.
\]
In the Fourier domain we conclude by the compact support of $\psi$
\begin{equation}
|\F\phi(u)|\lesssim\frac{1}{(1+|u|)^{\gamma}},\quad|\F\psi(u)|\lesssim\frac{|u|^{\gamma}}{(1+|u|)^{2\gamma}},\quad u\in\R.\label{eq:fPsi}
\end{equation}

\begin{lmm}
\label{lem:trick}Grant Assumption~\ref{assDiff} and let the
compactly supported father and mother wavelet functions $\phi$ and $\psi$ satisfy (\ref{eq:fPsi})
for some $\gamma>1$. Then for any $m\in\Z$, $j<l$ and $k\in K_j$
\begin{align*}
|\langle\psi_{l,m},\psi_{j,k}'\mu^{-1}\rangle| & \lesssim2^{l-(l-j)(\gamma-1/2)}+2^{-l+j},\\
\sum_{m\in\Z}|\langle\psi_{l,m},\psi_{j,k}'\mu^{-1}\rangle| &
\lesssim2^{l-(l-j)(\gamma-3/2)}+1\quad\text{and}\\
\sum_{m\in\Z}|\langle\psi_{l,m},\psi_{j,k}'\mu^{-1}\rangle|^{2}
&\lesssim2^{2l-(l-j)(2\gamma-2)}+2^{-(l-j)},
\end{align*}
where we have to replace $j$ by $j_0$ on the right-hand side for $\psi_{-1,k}=\phi_{j_0,k}$.
\end{lmm}
\begin{proof}
Let $\Gamma>0$ be large enough such that $\supp\phi\cup\supp\psi\subset[-\Gamma,\Gamma]$.
Noting that the following scalar product can only be nonzero if the support of $\psi_{l,m}$ is contained in $D$, a Taylor expansion of $\mu^{-1}$ yields for $j\ge j_0$
\begin{align*}
&\quad|\langle\psi_{l,m}\mu^{-1},\psi_{j,k}'\rangle|\\
&\le 
|\mu(2^{-l}m)|^{-1}|\langle\psi_{l,m},\psi'_{j,k}\rangle|+2^{-l}\Gamma\max_{
x:|x-m2^{-l}|\le2^{-l}\Gamma}|(\mu^{-1})'(x)|\int|\psi_{l,m}(x)||\psi{}_{j,k}'(x)|\d x\\
&\le 
\|\mu^{-1}\|_{L^\infty(D)}|\langle\psi_{l,m},\psi'_{j,k}\rangle|+2^{-l+j}\Gamma\max_{x:|x-m2^{-l}|\le2^{-l}\Gamma}|(\mu^{-1})'(x)|\|\psi\|_{L^{2}}\|\psi'\|_{L^{2}}.
\end{align*}
We conclude 
\begin{equation}
|\langle\psi_{l,m},\psi_{j,k}'\mu^{-1}\rangle|\lesssim\|\mu^{-1}\|_{L^\infty(D)}
|\langle\psi_{l,m},\psi'_{j,k}\rangle|+2^{-l+j}\|(\mu^{-1})'\|_{L^\infty(D)}.\label{eq:inLemTrick}
\end{equation}
Using Plancherel's identity,
$\F\psi_{l,m}(u)=\F[2^{l/2}\psi(2^{l}\bull-m)](u)=2^{-l/2}e^{imu2^{-l}}
\F\psi(2^ { -l}u)$ and (\ref{eq:fPsi}), we obtain
\begin{align*}
|\langle\psi_{l,m},\psi_{j,k}'\rangle|
&\le \frac{2^{-(j+l)/2}}{2\pi}\int|\F\psi(2^{-j}u)u\F\psi(2^{-l}u)|\d u\\
&\lesssim 2^{-(j+l)/2}\int\frac{2^{-(j+l)\gamma}|u|^{2\gamma+1}}{(1+2^{-j}|u|)^{2\gamma}(1+2^{-l}|u|)^{2\gamma}}\d u\\
&\le 2^{-(j+l)/2-(j+l)\gamma}\int\frac{|u|\d u}{2^{-2\gamma j}(1+2^{-l}|u|)^{2\gamma}}
=2^{l-(l-j)(\gamma-1/2)}\int\frac{|v|\d v}{(1+|v|)^{2\gamma}},
\end{align*}
where we have substituted $v=2^{-l}u$ in the last line. Due to $\gamma>1$, the integral in the last display is finite
so that combining this bound with (\ref{eq:inLemTrick}) yields the assertions,
noting that by the compact support of $\psi$ only for $\mathcal{O}(2^{l-j})$
many $m$ the scalar products $\langle\psi_{l,m},\psi_{j,k}'\mu^{-1}\rangle$
are nonzero.

For $j=-1$ we substitute again $v=2^{-l}u$ and obtain analogously
\begin{align*}
  |\langle\psi_{l,m},\psi'_{-1,k}\rangle|
  &\lesssim \int\frac{2^{-(j_0+l)/2-l\gamma}|u|^{\gamma+1}\d u} {(1+2^{-l}|u|)^{2\gamma}(1+2^{-j_0}|u|)^\gamma}
  \le2^{l-(l-j_0)(\gamma-1/2)}\int\frac{|v|\d v}{(1+|v|)^{2\gamma}}
\end{align*}
and 
\begin{align*}
  |\langle\psi_{l,m},\psi_{-1,k}'\mu^{-1}\rangle|
  &\lesssim\|\mu^{-1}\|_{L^\infty(D)}|\langle\psi_{l,m},\phi'_{j_0,k}\rangle|+2^{-l+j_0}\|(\mu^{-1})'\|_{L^\infty(D)}.\tag*{\qedhere}
\end{align*}
\end{proof}

Now we can bound the bias in (\ref{eq:DecompLin2}) in the multi-scale space $\mathcal{M}_{0}$.
\begin{lmm}\label{lem:biasMu}
  Let the weighting sequence $w$ be admissible and grant Assumption~\ref{assDiff} and (\ref{eq:fPsi})
  for some $\gamma\ge3/2+\delta,\delta\in(1,2]$. Then we have
  \[
  \big\|\big(\big\langle(\Id-\pi_{J+U})\mu,\psi_{j,k}'/(2\mu)\big\rangle\big)_{j\le J,k}\big\|_{\mathcal{M}}\lesssim2^{-J(s+1/2)}2^{-U(s+3/2)}w_{J}^{-1}.
  \]
\end{lmm}
\begin{proof}
Recall that we have by definition
\[
\big\|\big(\big\langle(\Id-\pi_{J^{+}})\mu,\psi_{j,k}'/(2\mu)\big\rangle\big)_{j\le J,k}\big\|_{\mathcal{M}}=\sup_{j\le J}\max_{k\in K_j}w_{j}^{-1}\big|\big\langle(\Id-\pi_{J^{+}})\mu,\psi_{j,k}'/(2\mu)\big\rangle\big|.
\]
As in the proof of Theorem~\ref{thrm:ucltMu} we have
$$\sup_{l,m:\supp \psi_{l,m}\cap S\neq\emptyset}2^{j(s+3/2)}|\langle\psi_{l,m},\mu\rangle|\lesssim \|\mu\|_{C^{s+1}(D)}.$$ 
Hence, for all $j\le J$,
\begin{align*}
\big|\big\langle(\Id-\pi_{J^{+}})\mu,\psi_{j,k}'/(2\mu)\big\rangle\big|&= 
\Big|\sum_{l>J^{+},m}\langle\mu,\psi_{l,m}\rangle\langle\psi_{l,m},\psi_{j,k}
'/(2\mu)\rangle\Big|\\
&\le 
\sup_{l>J^{+},m}2^{l}|\langle\mu,\psi_{l,m}\rangle|\sum_{l>J^{+},m}2^{-l}
|\langle\psi_{l,m},\psi_{j,k}'/(2\mu)\rangle|\\
&\lesssim 
2^{-J^{+}(s+1/2)}\sum_{l>J^{+},m}2^{-l}|\langle\psi_{l,m},\psi_{j,k}
'/(2\mu)\rangle|.
\end{align*}
Now Lemma~\ref{lem:trick} yields
\[
\sum_{l>J^{+},m}2^{-l}|\langle\psi_{l,m},\psi_{j,k}'/(2\mu)\rangle|\lesssim\sum_
{l>J^{+}}2^{-(l-j)(\gamma-3/2)}+\sum_{l>J^{+}}2^{-l}
\lesssim 2^{-(J^+-j)(\gamma-3/2)}+2^{-J^{+}}.
\]
Due to the monotonicity of $j\mapsto2^{j\delta}w_{j}^{-1}$, we conclude
for $\gamma\ge3/2+\delta$ 
\begin{align*}
\sup_{j\le J}\max_{k\in K_j}w_{j}^{-1}\big|\big\langle(\Id-\pi_{J^{+}})\mu,\psi_{j,k}
'/(2\mu)\big\rangle\big|&\lesssim  2^{-J^{+}(s+1/2)}\sup_{j\le
J}w_{j}^{-1}(2^{-(J^{+}-j)(\gamma-3/2)}+2^{-J^{+}})\\
&\lesssim  2^{-J^+(s+1/2)}2^{-(J^{+}-J)}w_{J}^{-1}.\tag*{\qedhere}
\end{align*}
\end{proof}
The second term in (\ref{eq:DecompLin2}) can be bounded by the following lemma.
\begin{lmm}\label{lem:ProjectionRest}
  Let the weighting sequence $w$ satisfy
  $\sqrt{j}2^{j}/w_{j}=\mathcal{O}(1)$ and grant Assumption~\ref{assDiff}
  and~(\ref{eq:fPsi}) for some $\gamma\ge5/2$. If $J_{n}^{+}=J_{n}+U_{n}$
  satisfies for some $\tau\in(0,1)$ 
  \[
  (\log n)^{2/\tau}n^{-1}2^{J_{n}^{+}}J_{n}^{+}=\mathcal{O}(1)\quad\text{and}\quad U_{n}\to\infty,
  \]
  then we have
  \[
  \Big\|\pi_{J_{n}}\Big(\big(\big\langle\mu_{n}-\mu,(\Id-\pi_{J_{n}^{+}})\Big(\frac{\psi_{j,k}'}{2\mu}\Big)\big\rangle\big)_{j,k}\Big)\Big\|_{\mathcal{M}}=o_{P}(n^{-1/2}).
  \]
\end{lmm}
\begin{proof}
In order to apply Proposition~\ref{prop:Concentration}, we need to calculate
the $L^{2}$-norm and the $L^{\infty}$-norm of $(\Id-\pi_{J_{n}^{+}})(\psi_{j,k}'/\mu)$.
For $j\in\{j_{0},\dots,J_n\}$ Parseval's identity and Lemma~\ref{lem:trick}
yield
\begin{align*}
\Big\|(\Id-\pi_{J_{n}^{+}})\Big(\frac{\psi_{j,k}'}{\mu}\Big)\Big\|_{L^{2}}^{2}&=
\sum_{l>J_{n}^{+},m}|\langle\psi_{l,m},\psi_{j,k}'/\mu\rangle|^{2}\\
&\lesssim 
\sum_{l>J_{n}^{+}}2^{2l-(l-j)(2\gamma-2)}+\sum_{l>J_{n}^{+}}2^{-2(l-j)+(l-j)_+}
\\
&\lesssim 
2^{j(2\gamma-2)}\sum_{l>J_{n}^{+}}2^{-l(2\gamma-4)}+2^{-(J_{n}^{+}-j)}\\
&\lesssim 
2^{-(J_{n}^{+}-j)(2\gamma-4)+2j}+2^{-(J_{n}^{+}-j)}\lesssim2^{-(J_{n}^{+}-j)}2^{
2j}.
\end{align*}
Another application of Lemma~\ref{lem:trick} yields
\begin{align*}
\Big\|(\Id-\pi_{J_{n}^{+}})\Big(\frac{\psi_{j,k}'}{2\mu}\Big)\Big\|_{\infty}&\le 
 \sum_{l>J_{n}^{+}}2^{l/2}\max_{n}|\langle\psi_{l,m},\psi_{j,k}'/\mu\rangle|\\
&\lesssim  
\sum_{l>J_{n}^{+}}2^{3l/2-(l-j)(\gamma-1/2)}+\sum_{l>J_{n}^{+}}2^{-l/2+j}\\
&\lesssim  
2^{-(J_{n}^{+}-j)(\gamma-2)+3j/2}+2^{-(J_{n}^{+}-j)/2+j/2}\lesssim2^{-(J_{n}^{+}
-j)/2}2^{3j/2}.
\end{align*}
The concentration inequality, Proposition~\ref{prop:Concentration}, yields for positive constants $c_{i}>0,i=1,2,\dots$,
\begin{align*}
 &\quad P\Big(\sup_{1\le j\le 
J_{n}}\max_{k}\sqrt{n}w_{j}^{-1}\big\langle\mu_{n}-\mu,(\Id-\pi_{J_{n}^{+}}
)\Big(\frac{\psi_{j,k}'}{2\mu}\Big)\big\rangle>t\Big)\\
&\le  \sum_{1\le j\le 
J_{n},k}P\Big(\sqrt{n}\big\langle\mu_{n}-\mu,(\Id-\pi_{J_{n}^{+}})\Big(\frac{
\psi_{j,k}'}{2\mu}\Big)\big\rangle>tw_{j}\Big)\\
&\lesssim  
\sum_{j=1}^{J_{n}}2^{j}\bigg(\exp(-c_{1}(2^{(J_{n}^{+}-j)/2}2^{-3j/2}\sqrt{n}w_{j
}t)^{\tau})\\
 & \qquad+\exp\Big(-\frac{c_{2}2^{(J_{n}^{+}-j)/2}t^{2}}{(2^{j}/w_{j})^{2}+t\max(2^{3j/2}(\log n)^{1/\tau},2^{j})/(w_{j}\sqrt{n})}\Big)\bigg).
\end{align*}
Since $jw_{j}^{-1}2^{3j/2}(\log n)^{1/\tau}n^{-1/2}\lesssim2^{j/2}j^{1/2}(\log n)^{1/\tau}n^{-1/2}\lesssim1$
and $J_{n}^{+}\lesssim\log n$ by the assumptions on $w_{j}$ and
$J_{n}^{+}$, we conclude for any $t>0$ and $n$ sufficiently large
\begin{align*}
 & \quad P\Big(\sup_{1\le j\le 
J_{n}}\max_{k}\sqrt{n}w_{j}^{-1}\big\langle\mu_{n}-\mu,(\Id-\pi_{J_{n}^{+}}
)\Big(\frac{\psi_{j,k}'}{2\mu}\Big)\big\rangle>t\Big)\\
&\lesssim  \sum_{j=1}^{J_{n}}2^{j}\Big(\exp(-c_{3}(\log
n)2^{\tau(J_{n}^{+}-j)/2}j^{\tau}t{}^{\tau})+\exp\Big(-c_{4}\frac{2^{(J_{n}
^{+}-j)/2}jt^{2}}{1+t}\Big)\Big)\\
&\lesssim  
\sum_{j=1}^{J_{n}}\exp\Big(j\big(\log2-c_{3}(J_{n}^{+})^{\tau-1}t^{\tau}\log 
n\big)\Big)+\sum_{j=1}^{J_{n}}\exp\Big(j\big(\log2-c_{5}2^{U_{n}/2}\frac{t^{2}}{
1+t}\big)\Big)\\
&\lesssim  e^{\log2-c_{3}(J_{n}^{+})^{\tau-1}t^{\tau}\log 
n}\frac{1-e^{J_{n}\big(\log2-c_{3}(J_{n}^{+})^{\tau-1}t^{\tau}\log 
n\big)}}{1-e^{\log2-c_{3}(J_{n}^{+})^{\tau-1}t^{\tau}\log n}}\\
 & \quad+e^{\log2-c_{6}2^{U_{n}/2}(t^{2}\wedge t)}\frac{1-e^{J_{n}(\log2-c_{6}2^{U_{n}/2}(t^{2}\wedge t))}}{1-e^{\log2-c_{6}2^{U_{n}/2}(t^{2}\wedge t)}}\\
&\lesssim  e^{-c_{3}(J_{n}^{+})^{\tau-1}t^{\tau}\log 
n}+e^{-c_{6}2^{U_{n}/2}(t^{2}\wedge t)}\to0.
\end{align*}
Finally note that all bounds hold true for the scaling function $\phi_{j_0,\cdot}$ if $j$ is replace by $j_0$.
\end{proof}

Now we have all pieces together to prove the multi-scale central
limit theorem.
\begin{proof}[Proof of Theorem~\ref{thrm:ucltDrift}]
Since $b$ has H\"older regularity $s>0$ on $S\subset D$, the bias can be 
bounded
by
\[
\|b-\pi_{J_{n}}(b)\|_{\mathcal{M}}=\sup_{j>J_{n}}\max_{k\in K_j}w_{j}^{-1}
|\langle\psi_{j,k},b\rangle|\lesssim\sup_{j>J_{n}}w_{j}^{-1}2^{-j(s+1/2)}=o(n^{
-1/2}).
\]
Using that the $\mathcal{M}_{0}$-norm is weaker than the $L^{\infty}(S)$-norm,
Lemma~\ref{lem:linearisation} and decomposition (\ref{eq:DecompLin2})
together with Lemmas~\ref{lem:biasMu} and \ref{lem:ProjectionRest}
yield 
\[
\big(\langle\hat{b}_{J_{n},U_{n}}-b,\psi_{j,k}\rangle\big)_{j\le
J_{n},k}=-\big(\big\langle\mu_{n}-\mu,\frac{\psi_{j,k}'}{2\mu}\big\rangle\big)_{
j\le
J_{n},k}+\tilde{R}_{J_{n},U_{n}}\quad\text{with }\|\pi_{J_{n}}(\tilde{R}_{
J_ { n },U_{n}})\|_{\mathcal{M}}=o_{P}(n^{-1/2}).
\]
Therefore, it remains to show that
\[
\beta_{\mathcal{M}_{0}}\big(\mathcal{L}\big(\pi_{J_{n}}\big((-\sqrt{n}
\langle\mu_ { n }
-\mu,\psi_{j,k}'/(2\mu)\rangle)_{j,k}\big)\big),\mathcal{L}(\mathbb{G}_{b}
)\big)\to0.
\]
This follows exactly as in Theorem~\ref{thrm:ucltMu}, where we use
that the factor $2^{j}$, by which the norms
\begin{equation}
\|\psi_{j,k}'/(2\mu)\|_{L^{2}}\lesssim2^{j},\quad\|\psi_{j,k}'/(2\mu)\|_{\infty}\lesssim2^{3j/2}\label{eq:bounds}
\end{equation}
are larger than $\|\psi_{j,k}\|_{L^{2}}$ and $\|\psi_{j,k}\|_{\infty}$,
respectively, is counterbalanced through the additional growth of the admissible weighting sequence $w$.
\end{proof}

\section{Remaining proofs}

\subsection{Proof of Theorem~\ref{thrm:deltaMethod}}\label{sec:proofDeltaMethod}

  \emph{Step 1:} For $\delta\in(0,1/2)$ and $0<c<C<\infty$ define \[\mathbb V_\xi:=\{\mu\in C^{7/4+\delta/2}(D):0<c<\mu\text{ and
}\|\mu\|_{C^{7/4+\delta/2}}<C\},\] $\mathbb V:=\mathcal M_0(w)$ and $\mathbb W:=\mathcal M_0(\tilde w)$ with
$w_j\le 2^{\delta j}$ and $\tilde w\ge 2^{(1+\delta)j}$. We first establish the Hadamard differentiability of
\begin{align*}
 \xi:\mathbb V_\xi  \subset \mathbb V \to \mathbb W,
\mu\mapsto \frac{\mu'}{2\mu}
\end{align*}
with derivative given by \eqref{eq:hadamardDerivative}.

To this end let $h\in\mathcal M_0(w)$ and $h_t\to h$ as $t\to0$.
For all $h_t$ such that $\mu+th_t$ is contained in~$\mathbb V_\xi$ for small
$t>0$ we have
\begin{align*}
&\quad \left\|\frac{\xi(\mu+th_t)-\xi(\mu)}{t}
-\xi_\mu'(h)\right\|_{\mathcal M(\tilde w)}
=\left\|\frac{\mu(\mu'+th_t')-\mu'(\mu+th_t)}{2(\mu+th_t)\mu t}-\frac{\mu
h'-\mu' h}{2\mu^2}
\right\|_{\mathcal M(\tilde w)}\\
&=\left\|\frac{\mu h_t'-\mu' h_t}{2(\mu+th_t)\mu }-\frac{\mu
h'-\mu' h}{2\mu^2}
\right\|_{\mathcal M(\tilde w)}
=\left\|\frac{\mu^2(h_t'-h')+\mu\mu'(h-h_t)+th_t(\mu'h-\mu
h')}{2\mu^2(\mu+th_t)}
\right\|_{\mathcal M(\tilde w)},\\
\intertext{using $\tilde w_j\ge 2^{j(1+\delta)}$ for
$\delta\in(0,1/2)$, this is bounded by}
&\left\|\frac{(h_t'-h')}{2(\mu+th_t)}
\right\|_{B_{\infty \infty}^{-3/2-\delta}}
+\left\|\frac{\mu'(h-h_t)}{2\mu(\mu+th_t)}
\right\|_{B_{\infty \infty}^{-3/2-\delta}}
+\left\|\frac{
th_t(\mu'h-\mu
h')}{2\mu^2(\mu+th_t)}
\right\|_{\mathcal M(\tilde w)}.\\
\intertext{Applying a pointwise multiplier theorem
\cite[Thm. 2.8.2]{triebel1983} and the continuous embedding 
$C^\rho\to B^\rho_{\infty \infty}$, we obtain up to constants the upper bound}
&\left\|\frac{1}{2(\mu+th_t)}
\right\|_{C^{7/4+\delta/2}}
\left\|h_t'-h'\right\|_{B_{\infty \infty}^{-3/2-\delta}}
+\left\|\frac{\mu'}{2\mu(\mu+th_t)}\right\|_{C^{7/4+\delta/2}}
\left\|h-h_t\right\|_{B_{\infty \infty}^{-1/2-\delta}}\\
&\quad+t\left\|\frac{
h_t(\mu'h-\mu
h')}{2\mu^2(\mu+th_t)}
\right\|_{\mathcal M(\tilde w)}\\
&\lesssim 
\left\|h_t-h\right\|_{B_{\infty \infty}^{-1/2-\delta}}
+\left\|h-h_t
\right\|_{B_{\infty \infty}^{-1/2-\delta}}+t
\lesssim \left\|h_t-h\right\|_{\mathcal M(w)}+t,
\end{align*}
where we have used $w_j\le 2^{j\delta}$ in the last step.
The last expression tends to 0 as $t\to0$ and this shows the Hadamard
differentiability of
$\xi\colon\mathbb V_\xi\to\mathbb W$.

  \emph{Step 2:}
To apply the delta method it is now important that $\hat\mu_{J_n}$ maps into $\mathbb
V_\xi$.
Theorem~\ref{thrm:ucltMu} gives conditions such that
$\|\hat\mu_{J_n}-\mu\|_{\mathcal
M(w)}=O(n^{-1/2})$. Provided that $2^{(9/4+\delta/2)J_n}w_{J_n}n^{-1/2}=o(1)$
we deduce from the fact that
$\hat\mu_{J_n}$ is developed until level $J_n$ only and from the ratio of the weights
at level~$J_n$ that $\|\hat\mu_{J_n}-\mu\|_{C^{7/4+\delta/2}}=o(1)$.
We conclude that with probability tending to one $\hat \mu\in\mathbb
V_\xi$. By modifying $\hat\mu$ on events with probability tending to zero we
can achieve that always $\hat\mu\in\mathbb V_\xi$. On the above assumptions
 we obtain the weak convergence $\sqrt{n}(\hat
\mu_{J_n}-\mu)\to\mathbb{G}_{\mu}$ in $\mathcal{M}_{0}(w)$ by
Theorem~\ref{thrm:ucltMu} and application of the delta method yields
the assertion.\qed

\subsection{Proof of Lemma~\ref{lem:Lepski}\label{sec:ProofofLepski}}

We will prove that:
\begin{enumerate}
\item for any $\tau\in(0,1)$ there are constants $0<c,C<\infty$ depending
only on $\tau,K,\psi$ such that for any $J\in\mathcal{J}_{n}$ satisfying
$J>J_{n}^{*}$ and for all $n\in\N$ large enough
\[
P(\hat{J}_{n}=J)\le C(n^{-cJ^{\tau}}+e^{-cJ}),
\]
\item there is an integer $M>0$ depending only on $d_{1},d_{2},K$ and
constants $0<c',C'<\infty$ depending on $\tau,K,\psi$ such that for
any $J\in\mathcal{J}_{n}$ satisfying $J<J_{n}^{*}-M$ and for all $n\in\N$ large enough
\[
P(\hat{J}_{n}=J)\le C'\big(n^{-c'J^{\tau}}+e^{-c'J}\big).
\]
\end{enumerate}
Given (i) and (ii), we obtain, for a constant $c''>0$, 
\begin{align*}
P\big(\hat{J}_{n}\notin[J_{n}^{*}-M,J_{n}^{*}]\big)&\le  
\sum_{J=J_{min}}^{J_{n}^{*}-M-1}P(\hat{J}_{n}=J)+\sum_{J=J_{n}^{*}+1}^{J_{max}}
P(\hat{J}_{n}=J)\\
&=  
\mathcal{O}\Big((J_{max}-J_{min})\big(n^{-c''(J_{min})^{\tau}}+e^{-c''J_{min}}
\big)\Big)\\
&=  \mathcal{O}\Big(\log 
n\big(n^{-c''J_{min}^{\tau}}+e^{-c''J_{min}}\big)\Big).
\end{align*}
Since $n^{-c''J_{min}^{\tau}}+e^{-c''J_{min}}$ decays polynomially
in $n$, the assertion of the lemma follows.

To show (i) and (ii), recall $J^{+}=J+U_{n}=J+\log_{2}\log n$. For
notational convenience we define 
\[
V(n,j):=(2^{3j}j/n)^{1/2},
\]
which is the order of magnitude of the stochastic error for projection level $j$.
Recall that for any $f\in\mathcal{M}\cap V_{J}$ we can bound 
\[
\|f\|_{L^\infty([a,b])}\lesssim\sum_{j\le J}2^{j/2}\max_{k\in K_j}|\langle f,\psi_{j,k}\rangle|\le\|f\|_{\mathcal{M}}\sum_{j\le J}2^{j/2}w_{j}\lesssim\|f\|_{\mathcal{M}}2^{J/2}w_{J}.
\]
Since any $J\in\mathcal{J}_{n}$ satisfies $2^{-J^{+}(s+1)}=o(1)$
and $(\log n)^{2/\tau}n^{-1}2^{J^{+}}J^{+}=\mathcal{O}(1)$, we conclude
from decomposition (\ref{eq:DecompLin2}) as well as from Lemma~\ref{lem:linearisation}
and Lemma~\ref{lem:ProjectionRest} (applied to $w_{j}=j^{1/2}2^{j}$) that
\[
\hat{b}_{J}-b=-\sum_{j\le J,k\in\Z}\Big\langle\mu_{n}-\mu,\frac{\psi_{j,k}'}{2\mu}\Big\rangle\psi_{j,k}+\underbrace{\sum_{j\le J,k\in\Z}\Big\langle(\Id-\pi_{J^+})\mu,\frac{\psi_{j,k}'}{2\mu}\Big\rangle\psi_{j,k}}_{=:B^{\mu}(J)}+\underbrace{\sum_{j>J,k\in\Z}\langle b,\psi_{j,k}\rangle\psi_{j,k}}_{=:B^{b}(J)}+\bar{R}_{J}
\]
for some remainder $\bar{R}_{J}\in V_{J}$ with 
\[
\|\bar{R}_{J}\|_{L^\infty([a,b])}=\mathcal{O}_{P}\big(n^{-1}J^{+}2^{2J^{+}}\big)+\mathcal{O}\big(2^{-J^{+}(2s+1)}\big)+o_{P}\big(n^{-1/2}J^{1/2}2^{-3J/2}\big).
\]
Moreover, Lemma~\ref{lem:biasMu} and Assumption~\ref{ass:SelfSim}
yield
\begin{equation}
\|B^{\mu}(J)+B^{b}(J)\|_{L^\infty([a,b])}\le(d_{2}+o(1))2^{-Js}.\label{eq:biasEst}
\end{equation}
Using $n^{-1/2}2^{2J^{+}}J^{+}=o(1)$ for all $J\in\mathcal{J}_{n}$,
we conclude
\begin{equation}
\|\hat{b}_{J}-b\|_{L^\infty([a,b])}\le\Big\|\sum_{j\le J,k\in K_j}\Big\langle\mu_{n}-\mu,\frac{\psi_{j,k}'}{2\mu}\Big\rangle\psi_{j,k}\Big\|_{L^\infty([a,b])}+(d_{2}+o(1))2^{-Js}+o_{P}(V(n,J)).\label{eq:ErrEst}
\end{equation}

With this preparation at hand we can proceed similarly as in \cite[Lem.
2]{gineNickl2010}. 

\emph{Part (i):} For any fixed $J>J_{n}^{*}$ we have 
\[
P\big(\hat{J}_{n}=J\big)\le\sum_{L\in\mathcal{J}_{n},L\ge J}P\big(\|\hat{b}_{J-1}-\hat{b}_{L}\|_{L^\infty([a,b])}>K\,V(n,L)\big).
\]
As in the derivation of (\ref{eq:ErrEst}) we obtain for $n$ sufficiently large 
\begin{align*}
\|\hat{b}_{J-1}-\hat{b}_{L}\|_{L^\infty([a,b])}&\le  \Big\|\sum_{J<j\le 
L,k}\Big\langle\mu_{n}-\mu,\frac{\psi_{j,k}'}{2\mu}\Big\rangle\psi_{j,k}\Big\|_{
L^\infty([a,b])}\\
 & \quad+\big(d_{2}+1\big)(2^{-(J-1)s}+2^{-Ls})+\frac{1}{4}(V(n,J-1)+V(n,L)).
\end{align*}
By definition of $J_{n}^{*}$ we have for any $L>J>J_{n}^{*}$ that
\begin{align*}
\big(d_{2}+1\big)(2^{-(J-1)s}+2^{-Ls})&\le  2(d_{2}+1)2^{-J_{n}^{*}s}
\le  \frac{K}{2}V(n,J_{n}^{*})\le\frac{K}{2}V(n,L).
\end{align*}
Therefore,
\[
P\big(\hat{J}_{n}=J\big)\le\sum_{L\in\mathcal{J}_{n},L\ge J}P\Big(\Big\|\sum_{J<j\le L,k}\Big\langle\mu_{n}-\mu,\frac{\psi_{j,k}'}{2\mu}\Big\rangle\psi_{j,k}\Big\|_{L^\infty([a,b])}>\frac{K-1}{2}V(n,L)\Big).
\]
Analogously to (\ref{eq:ApplyBernstein}) and using (\ref{eq:bounds}),
Proposition~\ref{prop:Concentration} yields for any $\tau\in(0,1)$
and constants $c_{1},...,c_{4}>0$:
\begin{align}
 &\quad P\Big(\Big\|\sum_{J<j\le L,k\in 
K_j}\Big\langle\mu_{n}-\mu,\frac{\psi_{j,k}'}{2\mu}\Big\rangle\psi_{j,k}\Big\|_{
L^\infty([a,b])}>\frac{K-1}{2}V(n,L)\Big)\nonumber \\
&\le  P\Big(\sum_{J<j\le
L}2^{j/2}\max_{k\in K_j}\Big|\Big\langle\mu_{n}-\mu,\frac{\psi_{j,k}'}{2\mu}
\Big\rangle\Big|>\frac{K-1}{2}\sqrt{\frac{2^{3L}L}{n}}\Big),\nonumber \\
\intertext{using that $\sum_{k=0}^\infty2^{-k/2}\le7/2$, we obtain the upper
bound}
& \quad P\Big(\max_{J<j\le
L,k}n^{1/2}2^{-L}\Big|\Big\langle\mu_{n}-\mu,\frac{\psi_{j,k}'}{2\mu}
\Big\rangle\Big|>\frac{K-1}{7}L^{1/2}\Big)\nonumber \\
&\le  \sum_{J<j\le 
L,k}P\Big(n^{1/2}2^{-j}\Big|\Big\langle\mu_{n}-\mu,\frac{\psi_{j,k}'}{2\mu}
\Big\rangle\Big|>\frac{K-1}{7}j^{1/2}\Big)\nonumber \\
&\le  c_{1}\sum_{J<j\le L}\big(e^{-c_{2}j^{\tau}\log n}+e^{-c_{3}j}\big)
\le c_{1}(n^{-c_{4}J^{\tau}}+e^{-c_{4}J}),\label{eq:prLemLep}
\end{align}
where we require that $K$ is chosen sufficiently large, depending on $\|\psi'_{j,k}/\mu\|_{L^\infty(S)}$ and $\Sigma_{j,k}$.
It remains to sum this upper bound over all $L\in\mathcal{J}_{n}$
with $L\ge J$ which yields the claim since $\mathcal{J}_{n}$ contains
no more than $\log n$ elements.

\emph{Part (ii):} Let $J<J_{n}^{*}-M$ for some $M\in\N$ to be specified
below. We have
\[
P\big(\hat{J}_{n}=J\big)\le P\big(\|\hat{b}_{J}-\hat{b}_{J_{n}^{*}}\|_{L^\infty([a,b])}\le K\,V(n,J_{n}^{*})\big).
\]
Using Assumption~\ref{ass:SelfSim} and the triangle inequality,
we obtain similarly to (\ref{eq:ErrEst}), for sufficiently large~$n$,
\begin{align*}
\|\hat{b}_{J}-\hat{b}_{J_{n}^{*}}\|_{L^\infty([a,b])}&\ge  
\frac{d_{1}}{2}2^{-Js}-\big(d_{2}+1\big)2^{-J_{n}^{*}s}\\
 &\quad -\Big\|\sum_{J<j\le J_{n}^{*},k\in K_j}\Big\langle\mu_{n}-\mu,\frac{\psi_{j,k}'}{2\mu}\Big\rangle\psi_{j,k}\Big\|_{L^\infty([a,b])}-\frac{1}{4}(V(n,J_{n}^{*})+V(n,J)).
\end{align*}
Owing to $J<J_{n}^{*}-M$, $s\ge1$ and the definition of $J_{n}^{*}$,
we can bound
\begin{align*}
\frac{d_{1}}{2}2^{-Js}-\big(d_{2}+1\big)2^{-J_{n}^{*}s}&\ge  
(d_{2}+1)\Big(\frac{d_{1}}{2(d_{2}+1)}2^{M-1}-\frac{1}{2}\Big)2^{-(J_{n}^{*}
-1)s}\\
&\ge  
\frac{K}{4}\Big(\frac{d_{1}}{2(d_{2}+1)}2^{M-1}-\frac{1}{2}\Big)V(n,J_{n}^{*}
-1)\\
&\ge 
\frac{K}{16}\Big(\frac{d_{1}}{2(d_{2}+1)}2^{M-1}-\frac{1}{2}\Big)V(n,J_{n}^{*}
),
\end{align*}
where we have used in the last inequality that $(J_{n}^{*}-1)/J_{n}^{*}\ge1/2$
for $n$ sufficiently large. We conclude
\begin{align*}
\|\hat{b}_{J}-\hat{b}_{J_{n}^{*}}\|_{L^\infty([a,b])}
&\ge\tilde{K}\,V(n,J_{n}^{*})-\Big\|\sum_{J<j\le J_{n}^{*},k\in 
K_j}\Big\langle\mu_{n}-\mu,\frac{\psi_{j,k}'}{2\mu}\Big\rangle\psi_{j,k}\Big\|_{
L^\infty([a,b])}
\end{align*}
with $\tilde K:=\frac{Kd_{1}}{32(d_{2}+1)}2^{M-1}-\frac{K}{32}-\frac{1}{2}$.
Since $\tilde{K}>K$ for $M$ large enough, we obtain similarly as 
in~(\ref{eq:prLemLep}) for any $\tau\in(0,1)$ and some $c',C'>0$ 
\begin{align*}
P\big(\hat{J}_{n}=J\big) & \le P\Big(\Big\|\sum_{J<j\le J_{n}^{*},k}\Big\langle\mu_{n}-\mu,\frac{\psi_{j,k}'}{2\mu}\Big\rangle\psi_{j,k}\Big\|_{L^\infty([a,b])}\ge(\tilde{K}-K)V(n,J_{n}^{*})\Big)\\
 & \le C'\big(n^{-c'J^{\tau}}+e^{-c'J}\big).\tag*{\qed}
\end{align*}

\subsection{Proof of Corollary~\ref{cor:sHat}\label{sec:ProofSHat}}

Owing to $(cn/\log n)^{1/(2s+3)}\le2^{J_{n}^{*}}\le(Cn/\log n)^{1/(2s+3)}$
for constants $0<c<C$, we find 
\[
\frac{\log n-\log\log n}{2(\log2)J_{n}^{*}}+\frac{\log c}{2(\log2)J_{n}^{*}}-\frac{3}{2}\le s\le\frac{\log n-\log\log n}{2(\log2)J_{n}^{*}}+\frac{\log C}{2(\log2)J_{n}^{*}}-\frac{3}{2}.
\]
Since $P(\hat{J}_{n}\le J_{n}^{*})\to1$ by Lemma~\ref{lem:Lepski}
and due to $v_{n}^{-1}=o_P(1)$, we obtain with probability
converging to one that
\[
s\ge\max\bigg(1,\,\frac{\log n-\log\log n}{2(\log2)(\hat{J}_{n}+v_{n})}-o_P\Big(\frac{v_{n}}{\hat{J}_{n}}\Big)-\frac{3}{2}\bigg)\ge\hat{s}_{n}.
\]
Moreover, since $P(J_{n}^{*}-M\le\hat{J}_{n}\le J_{n}^{*})\to1$ and $v_{n}^{-1}=o_P(1)$,
we have with probability converging to one
\[
s-\hat{s}_{n}\le\frac{\log n-\log\log n}{2(\log2)J_{n}^{*}}\Big(1-\frac{1}{1+v_{n}/J_{n}^{*}}\Big)+\frac{\log C}{2(\log2)J_{n}^{*}}+\frac{3v_{n}}{2(J_{n}^{*}-M)}\lesssim\frac{v_{n}}{J_{n}^{*}}.\tag*{\qed}
\]

\bibliographystyle{apalike}
\bibliography{bib}

\end{document}